\documentclass[a4paper,10pt,onecolumn]{article}

\usepackage[tmargin=2cm,lmargin=2cm,rmargin=2cm,bmargin=2.5cm]{geometry}                		
\usepackage{cite}
\usepackage{amsmath,amssymb,amsfonts}
\usepackage{algorithmic}
\usepackage{graphicx}
\usepackage{textcomp}
\usepackage{optidef}
\usepackage{hyperref}

\usepackage{comment}
\usepackage{mathtools}
\usepackage{enumerate}
\usepackage{overpic}
\usepackage{xcolor}
\usepackage{todonotes}

\usepackage{siunitx}
\usepackage{makecell}
\usepackage{multirow}
\allowdisplaybreaks[4]
\usepackage{cuted}
\def\fuel{\text{fuel}}
\def\gt{\text{gt}}
\def\e{\text{el}}

\def\gen{\text{gen}}
\def\el{\text{em}}
\def\b{\text{b}}
\def\c{\text{c}}
\def\drv{\text{drv}}

\DeclareMathOperator{\proj}{\pi}
\DeclareMathOperator{\diag}{diag}
\newcommand{\vvec}[1]{\overrightarrow{#1}}

\newcommand\mdsedit[1]{\textcolor{black}{#1}}
\newcommand\mdseditbis[1]{\textcolor{black}{#1}}
\newcommand\mdsedittris[1]{\textcolor{black}{#1}}
\newcommand\finaledit[1]{\textcolor{black}{#1}}

\begin{document}
\title{Predictive energy management for hybrid electric aircraft propulsion systems}

\author{Martin Doff-Sotta\thanks{M.~Doff-Sotta, M.~Cannon and M.~Bacic (on part-time secondment from Rolls-Royce plc) are with the Control Group, University of Oxford, Parks Road, Oxford, OX1 3PJ, United Kingdom  (e-mail: martin.doff-sotta@eng.ox.ac.uk, mark.cannon@eng.ox.ac.uk, marko.bacic@eng.ox.ac.uk).},~~Mark Cannon,~~Marko Bacic}

\maketitle

\begin{abstract}
We present a Model Predictive Control (MPC) algorithm for energy management in aircraft with hybrid electric propulsion systems consisting of gas turbine and electric motor components. Series and parallel configurations are considered.
By combining a point-mass aircraft dynamical model with models of electrical losses and losses in the gas turbine, the fuel consumed over a given future flight path is minimised subject to constraints on the battery, electric motor and gas turbine. The optimization is formulated as a convex problem under mild assumptions and its solution is used to define a predictive energy management control law that takes into account the variation in aircraft mass during flight.
%
We investigate the performance of algorithms for solving this problem. An Alternating Direction Method of Multipliers (ADMM) algorithm is proposed and compared with a general purpose convex interior point solver.  
We also show that the ADMM implementation reduces the required computation time by orders of magnitude in comparison with a general purpose nonlinear programming solver, making it suitable for real-time supervisory energy management control.



%
%
\noindent\textbf{keywords:}
Alternating Direction Method of Multipliers (ADMM), Convex Programming, Energy Management, Hybrid Aircraft, Model Predictive Control (MPC).
\end{abstract}

\section{Introduction}
\label{sec:introduction}
Aviation currently contributes to around 2\% of current world-wide human-made CO2 emissions, but demand for air travel is predicted to grow significantly.
The aviation industry is committed to realising this growth sustainably with a drastic reduction of CO2 emissions by 2050.
One avenue identified to achieve this ambitious goal is the development of greener aviation based on new propulsion concepts.  

Aircraft equipped with turbo-electric and hybrid electric powertrains are considered in~\cite{welstead2016conceptual, lammen2020energy, antcliff2017conceptual} where it is shown that reductions in emissions and energy savings can potentially be achieved. In~\cite{welstead2016conceptual}, simulations of a commercial airliner with boundary layer ingestion and a turbo-electric propulsion system predict mission fuel burn savings of up to 7\% relative to a conventional propulsion system. A distributed electric propulsion concept for the transonic cruise range proposed in~\cite{schmollgruber2020multidisciplinary} is likewise expected to provide a 7\% reduction in fuel. Potential energy savings were demonstrated for a concept year-2030 aircraft equipped with a parallel-hybrid propulsion system combined with an all-electric propulsion system in \cite{antcliff2017conceptual}. 

Hybrid-electric propulsion systems rely on energy management controllers to allocate power demand between the different components of the powertrain. The energy management problem can be tackled with heuristic strategies such as peak-shaving schemes \cite{strack2017conceptual}, charge-depleting-charge-sustaining policies \cite{doff2020optimal}, approaches based on state machines \cite{spagnolo2018finite} and rule-based fuzzy logic \cite{caux2010line}. More sophisticated suboptimal control strategies have also been proposed, for example using neural networks \cite{moreno2006energy} and neuro-fuzzy adaptive control \cite{karunarathne2010model}. 

Optimisation techniques that seek to minimise a cost function (such as fuel consumption) have also been proposed for energy management problems. For example, the so-called equivalent fuel consumption  minimisation strategy is widely used in hybrid fuel cell systems \cite{garcia2012viability}. 
Globally optimal policies have been computed offline using dynamic programming~\cite{bradley2009energy, kim2007power, lin03} but the  computation required is prohibitive for real-time implementation.
Other approaches based on $\mathcal{H}_\infty$ control~\cite{pisu2007comparative} and optimal adaptive control~\cite{lin2011energy} have also been proposed. A popular framework for energy management problems in electric and hybrid-electric ground vehicles is Model Predictive Control (MPC)~\cite{koot05, east2019energy, buerger2018fast}. The energy management problem is formulated as a  receding-horizon constrained optimisation problem, and an optimal power split is found at each discrete time step. Since MPC is a feedback control strategy that is updated with information on the system state at each time step, it can provide robustness to modelling uncertainty and prediction errors. 
Although MPC has been proposed for energy management problems in hybrid-electric aircraft \cite{seok2017coordinated, dunham2019distributed}, none of these approaches considered lossless convexification of the nonlinear programming problem.

A convex energy management formulation is proposed in~\cite{doff2020optimal}, which considers a parallel-hybrid aircraft with nonlinear constraints in a model predictive control framework.
Fuel consumption is predicted over a future flight profile and is minimised subject to constraints on state, trajectory and physical limitations of the components of the propulsion system. The associated receding-horizon nonlinear programming problem is posed as a convex program and solved using the general-purpose convex optimisation framework CVX \cite{cvx_software}.

This paper extends the results in~\cite{doff2020optimal} to series-hybrid architectures and describes a specialised ADMM solver for efficient online optimisation. The paper is organised as follows. Mathematical models of the powertrain components and aircraft dynamics are developed in Section~\ref{sec:modelling}. Energy management problems for series and parallel configurations are stated as receding-horizon optimisation problems in Section~\ref{sec:dt_model}.  Section~\ref{sec:convex} presents a series of simplifications that yield convex relaxations of these problems. In particular, a unified formulation is proposed for both powertrain configurations. The ADMM solver is presented in Section~\ref{sec:admm} and its performance and potential for real-time implementation are discussed in Section~\ref{sec:simulation_results}. Conclusions are presented in Section~\ref{sec:conclusion}.

\section{Modelling}\label{sec:modelling}

This section derives models of the aircraft dynamics and powertrain components (battery, electric motor, gas turbine etc.), which will be used to formulate the energy management problem as a model-based optimisation problem. 

We consider a hybrid electric aircraft propulsion system with either a series or parallel topology. When the power output demand is negative, which may occur for example while the aircraft is descending, we consider the possibility of using the same powertrain to generate electrical energy (i.e.\ operating in a ``windmilling'' mode) in order to recharge the battery. In practice a variable-pitch fan would be required for this functionality, which would increase complexity.  

\subsection{Aircraft Dynamics}

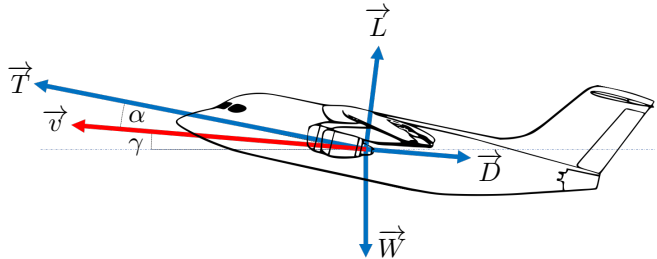
\begin{figure}[htb]
\centerline{\begin{overpic}[scale=0.35,trim=-1mm 0mm 0mm 8mm, clip]{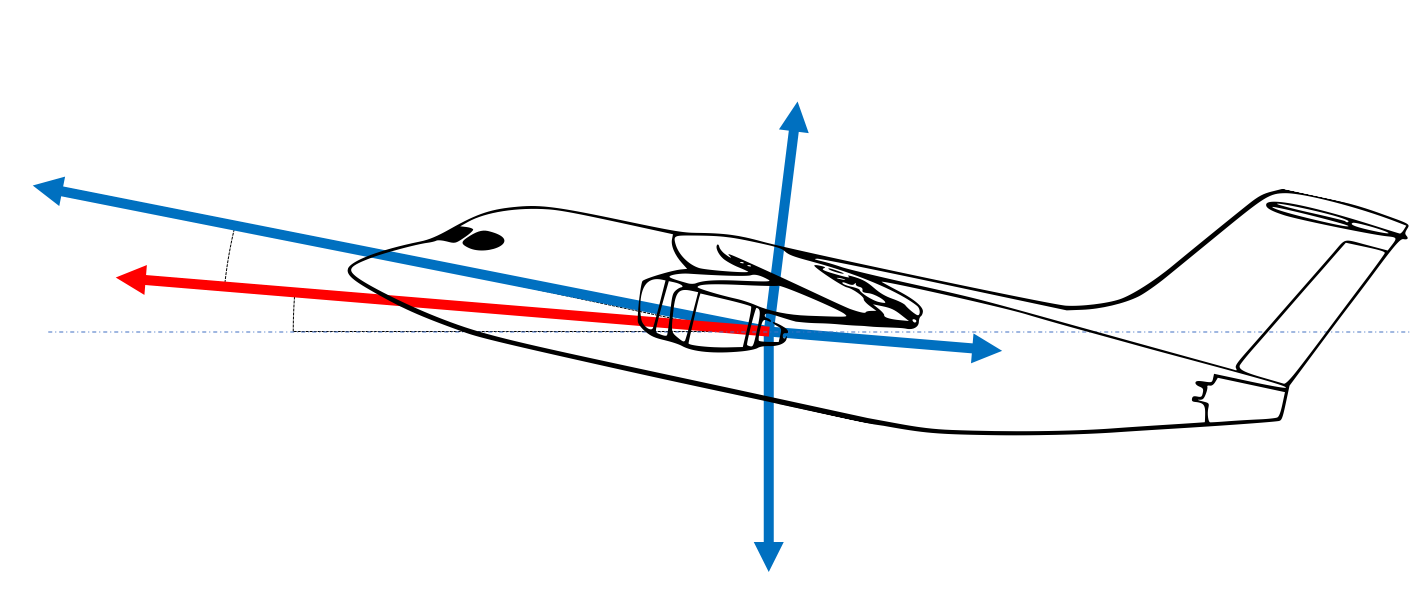}
\put(-1.5,28){$\vvec{T}$}
\put(72,14){$\vvec{D}$}
\put(56.5,2){$\vvec{W}$}
\put(54.5,37){$\vvec{L}$}
\put(4,22){$\vvec{v}$}
\put(17.5,23.3){$\alpha$}
\put(17.5,19){$\gamma$}
\end{overpic}}
\vspace{-2mm}
\caption{Aircraft forces and motion.}
\label{fig:aircraft}
\end{figure}

The aircraft motion is constrained by its dynamic equations. Assuming a point-mass model~\cite{stevens2015aircraft} and referring to Figure \ref{fig:aircraft}, the equilibrium of forces yields
\begin{equation*}
m\frac{\mathrm{d}}{\mathrm{d}t}( \vvec{v})  =\vvec{T} + \vvec{L} + \vvec{D} + \vvec{W}, 
\end{equation*}
where $ \vvec{v}$ is the velocity vector, $m$ the instantaneous mass of the aircraft, $\vvec{T}$ the vector of thrust, $\vvec{L}$ and $\vvec{D}$ are the lift and drag vectors and $\vvec{W}$ is the aircraft weight.

Using the coordinates ($v$,$\gamma$), where $v$ is the velocity vector magnitude and $\gamma$ is the flight path angle, and projecting the vector equation in wind axes along the drag vector $\vvec{D}$ yields
\[
 m\frac{\mathrm{d}}{\mathrm{d}t}v + m g \sin{\gamma}   = T \cos{\alpha} - \frac{1}{2}C_D\rho S v^2 .
\]
Here $S$ is the wing area, $\rho$ is the density of air, $g$ is acceleration due to gravity, $C_D=C_D(\alpha)$ the drag coefficient and $\alpha$ the angle of attack.
Projecting along the lift vector $\vvec{L}$ yields
\[
m v \frac{\mathrm{d}}{\mathrm{d}t} \gamma + mg \cos{\gamma}= T \sin{\alpha} + \frac{1}{2}C_L\rho S v^2 ,
\]
%
where $C_L=C_L(\alpha)$ is the lift coefficient. 

The drive power is given as follows
\[
P_{\drv}=\vvec{T} \cdot \vvec{v} = m\frac{\mathrm{d}}{\mathrm{d}t}(\frac{1}{2} v^2)  + \frac{1}{2}C_D\rho S v^3 +  m g v \sin{\gamma}  . 
\]
\subsection{Hybrid Propulsion System}

\subsubsection{Parallel architecture}
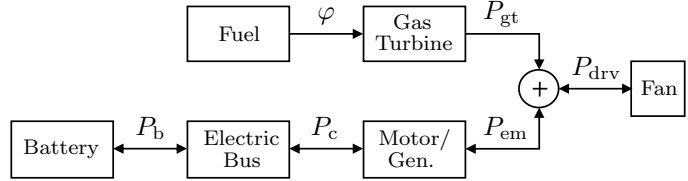
\begin{figure}[h!]
\centerline{\begin{overpic}[scale=0.644,trim=0mm 0mm 0mm -2mm,clip]{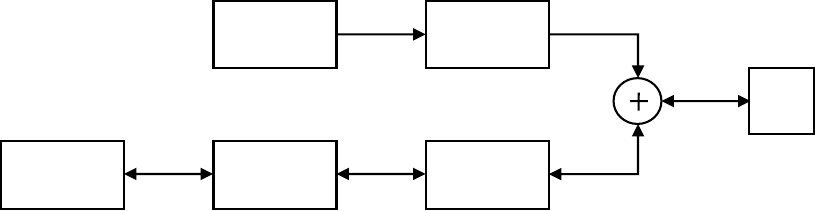}
\put(1.9,3.7){\footnotesize{Battery}}
\put(28.6,5){\footnotesize{Electric}}
\put(31.5,1.7){\footnotesize{Bus}}
\put(54.5,5){\footnotesize{Motor/}}
\put(55.9,1.7){\footnotesize{Gen.}}
\put(30.8,20.6){\footnotesize{Fuel}}
\put(56.7,22.3){\footnotesize{Gas}}
\put(54,19.2){\footnotesize{Turbine}}
\put(93.4,12.5){\footnotesize{Fan}}
\put(18.5,6.5){$P_{\b}$}
\put(44.5,6.5){$P_{\c}$}
\put(70,6.5){$P_{\el}$}
\put(45.5,23.8){$\varphi$}
\put(70,23.8){$P_{\gt}$}
\put(83.1,15.5){$P_{\drv}$}
\end{overpic}}
    \caption{Parallel-hybrid propulsion architecture.}
    \label{fig:propulsion}
\end{figure}

In the parallel architecture ($\mathcal{T} = \mathcal{P}$), a gas turbine producing power $P_{\gt}$ is mechanically coupled with an electric motor with power output $P_{\el}$ in a parallel arrangement  (Fig.~\ref{fig:propulsion}). 
These two power sources are combined to give the power output of the propulsion system, $P_{\drv}$, via
\[
P_{\drv}(t) = P_{\gt}(t) + P_{\el}(t), 
\]
where $100\%$ efficiency in drivetrain components is assumed.

\subsubsection{Series architecture} 
\begin{figure}[h!]
\centerline{\begin{overpic}[scale=0.23,trim=0mm 0mm 0mm -2mm,clip]{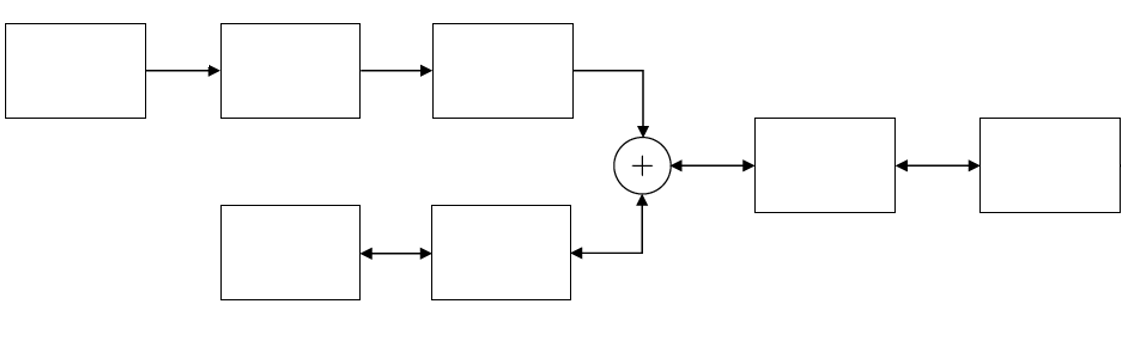}
\put(21,6.7){\footnotesize{Battery}}
\put(39.7,8.2){\footnotesize{Electric}}
\put(41.5,5){\footnotesize{Bus}}
\put(68.5,15.5){\footnotesize{Motor/}}
\put(69.5,12.5){\footnotesize{Gen.}}
\put(3.9,22.6){\footnotesize{Fuel}}
\put(23.5,24){\footnotesize{Gas}}
\put(21,21){\footnotesize{Turbine}}
\put(41.5,22.6){\footnotesize{Gen.}}
\put(90.4,14.5){\footnotesize{Fan}}
\put(32.3,9.5){$P_{\b}$}
\put(51,9.5){$P_{\c}$}
\put(60,16.7){$P_{\e}$}
\put(15.2,26){$\varphi$}
\put(32.3,26){$P_{\gt}$}
\put(51,26){$P_{\gen}$}
\put(79.1,16.7){$P_{\drv}$}
\end{overpic}}
\caption{Series-hybrid propulsion architecture.}
    \label{fig:propulsion_series}
\end{figure}
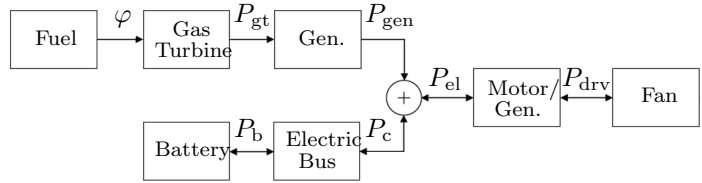

In the series architecture ($\mathcal{T} = \mathcal{S}$), the propulsion system power output $P_{\drv}$ is delivered by an electric motor taking electrical power $P_{\e}$ from two sources: a battery with effective power output $P_{\c}$ and a turbo generator set (gas turbine in series with an electric generator) with power output $P_{\gen}$ (Fig.~\ref{fig:propulsion_series}). The power balance is given by
\[
P_{\e}(t) = P_{\gen}(t) + P_{\c}(t). 
\]

\subsection{Battery}
The battery is modelled as an equivalent circuit with internal resistance $R$ and open-circuit voltage $U$, so that the input-output map between its chemical power $P_{\b}$ and the effectively delivered electrical power $P_{\c}$ is given by \cite{guzzella2007vehicle}
\begin{align*}
P_{\b}&=g\bigl(P_{\c}\bigr), \\
&=\frac{U^2}{2R}\left(1-\sqrt{1-\frac{4R}{U^2}P_{\c}}\right), 
\end{align*}
where $U$ and $R$ are assumed constant~\cite{east2018}. 
The evolution of the battery state of charge (SOC) $E(t)$ is given by
\begin{equation}
\dot{E}=-P_{\b}  
\label{eq:E}
\end{equation}
and $E(t)$ is subject at all times to upper and lower bounds 
\[
\underline{E} \leq E \leq  \overline{E}. 
\]
\subsection{Gas turbine}
The rate of change of mass of the aircraft is given by
\begin{equation}
 \dot{m}  = \varphi = -f(P_{\gt}(t), \omega_{\gt}(t)), 
\label{eq:fuel}
\end{equation}
where $\varphi$ is the rate of fuel consumption and $f(P_{\gt},\omega_{\gt})$ is a piecewise-quadratic function of the gas turbine power output $P_{\gt}$ and shaft rotation speed $\omega_{\gt}$. We assume that $f(\cdot)$ can be determined empirically from fuel map data in the form 
\begin{equation*}
\begin{split}
\varphi & = f(P_{\gt},\omega_{\gt}),  \\
& = \beta_{2}(\omega_{\gt})P_{\gt}^2 + \beta_{1}(\omega_{\gt})P_{\gt}  + \beta_{0}(\omega_{\gt}) , 
\end{split}
\end{equation*}
with $\beta_{2}(\omega_{\gt})\!\geq\! 0$ and $\beta_{1}(\omega_{\gt}) \!>\! 0$ in the operating range of $\omega_{\gt}$. 
The power $P_{\gt}$ and shaft rotation speed $\omega_{\gt}$ are limited by  
\begin{gather*}
\underline{P}_{\gt} \leq P_{\gt} \leq  \overline{P}_{\gt},
\\
\underline{\omega}_{\gt} \leq \omega_{\gt} \leq  \overline{\omega}_{\gt}. 
\end{gather*}
These limits apply to both parallel and series configurations, and in the latter case they constrain the turbo generator set.   
\subsection{Electric motor}
In the parallel configuration, the electric motor input-output map between input electrical power $P_{\c}$ and effective mechanical power $P_{\el}$ is modelled by a piecewise-quadratic  function 
\begin{equation*}
\begin{split}
P_{\c} &= h(P_{\el}(t), \omega_{\el}(t)), \\
& = \kappa_{2}(\omega_{\el})P_{\el}^2  + \kappa_{1}(\omega_{\el})P_{\el}  + \kappa_{0}(\omega_{\el}), 
\end{split}
\end{equation*}
where $\omega_{\el}$ is the electric motor shaft rotation speed and $\kappa_{2}(\omega_{\el}) \geq 0$, $\kappa_{1}(\omega_{\el}) > 0$ for all $\omega_{\el}$ in the operating range. The function $h(\cdot)$ can be determined empirically from electric motor loss data.
The limitations on the electric motor power and shaft rotation speeds are set by the following constraints
\begin{gather*}
\underline{P}_{\el} \leq P_{\el} \leq  \overline{P}_{\el},
\\
\underline{\omega}_{\el} \leq \omega_{\el} \leq  \overline{\omega}_{\el}. 
\end{gather*}

In the series configuration, the input-output map between the input electrical power $P_{\e}$ and effective mechanical power $P_{\drv}$ is likewise modelled by $P_{\e} = h(P_{\drv}(t), \omega_{\drv}(t))$,  where $\omega_{\drv}$ is the fan shaft rotation speed.
The limitations on the electric motor power and shaft rotation speeds are set by the following constraints
\begin{gather*}
\underline{P}_{\drv} \leq P_{\drv} \leq  \overline{P}_{\drv},
\\
\underline{\omega}_{\drv} \leq \omega_{\drv} \leq  \overline{\omega}_{\drv}. 
\end{gather*}

\subsection{Electric generator}
In the series configuration a generator converts the gas turbine mechanical power $P_{\gt}$ into electrical power $P_{\gen}$. This electrical machine is modelled by a piecewise-quadratic function
\begin{equation*}
\begin{split}
P_{\gt} &= h_{\gen}(P_{\gen}(t), \omega_{\gen}(t)), \\
& = \nu_{2}(\omega_{\gen})P_{\gen}^2  + \nu_{1}(\omega_{\gen})P_{\gen}  + \nu_{0}(\omega_{\gen}), 
\end{split}
\end{equation*}
where  $\omega_{\gen}$ is the electric generator shaft rotation speed and $\nu_{2}(\omega_{\gen}) \geq 0$, $\nu_{1}(\omega_{\gen}) > 0$ for all $\omega_{\gen}$ in the operating range. The loss map $h_{\gen}(\cdot)$ can be determined empirically from electric generator loss data. 
The limits on power and shaft rotation speed for the electric generator are encapsulated by the inequality constraints given for the gas turbine. 

\subsection{Objective}
The problem at hand is to find the real-time optimal power split between the gas turbine and electric motor that minimises
\begin{equation*}\label{eq:contin-time_objective}
    J = \int_{0}^{T}{f(P_{\gt}(t), \omega_{\gt}(t))}\mathrm{d}t,
\end{equation*}
while satisfying constraints on the battery SOC and limits on power flows throughout the powertrain, and while producing sufficient power to follow a prescribed flight path.



\section{Discrete-time optimal control}\label{sec:dt_model}

This section describes a discrete-time model that enables the optimal power split between battery and fuel over a given future flight path to be determined as a finite-dimensional optimisation problem.
For a fixed sampling interval $\delta$, we consider a predictive control strategy that minimises, at each sampling instant, the fuel consumption over the remaining flight path. 
The optimisation is performed subject to the dynamics of the aircraft mass and the battery SOC. The problem is also subject to limits on energy stored in the battery (to prevent deep discharging or overcharging) and limits on power flows corresponding to physical and safety constraints.

The optimal solution to the fuel minimisation problem at the $k$th sampling instant is computed using estimates of the battery SOC $E(k\delta)$ and the aircraft mass $m(k\delta)$, so that $E_0 = E(k\delta)$ and $m_0=m(k\delta)$ at any time $k\delta$. The control law at time $k\delta$ is defined by the first time step of this optimal solution.
The notation $\{ x_{0} , x_{1} , \ldots x_{N-1}\}$ is used for the sequence of current and future values of a variable $x$ predicted at the $k$th discrete-time step, so that $x_i$ denotes the predicted value of $x\bigl((k+i)\delta\bigr)$. The horizon $N$ is chosen so that $N = \lceil T/\delta \rceil - k$, and hence $N$ shrinks as $k$ increases and $k\delta$ approaches $T$.

The discrete-time approximation of the objective is  
\begin{equation}
    J = \sum_{i=0}^{N-1}{f_i(P_{\gt,i},\omega_{\gt, i})} \, \delta,
\label{eq:obj}
\end{equation}
with, for $i=0,\dots,N-1$,
\begin{align}
\mbox{}\hspace{-3mm} f_i(P_{\gt, i}, \omega_{\gt, i}) &= \beta_{2}(\omega_{\gt,i})P_{\gt,i}^2 \!+\! \beta_{1}(\omega_{\gt, i})P_{\gt, i}  \!+\! \beta_{0}(\omega_{\gt, i}), 
\label{eq:m_k0}
\\
m_{i+1} &= m_{i} - f_i(P_{\gt,i}, \omega_{\gt,i}) \, \delta , 
\label{eq:m_k}
\end{align}
where the forward \mdsedit{Euler approximation has been used to discretise~(\ref{eq:fuel}).
The same approach applied to (\ref{eq:E}) yields} the discrete-time battery model
\begin{align}
E_{i+1} &= E_{i} - P_{\b,i} \, \delta, 
\label{eq:E_k}
\\
 P_{\b,i} &= g_i (P_{\c,i}) , 
\nonumber 
\\
&= \frac{U^2}{2R}\biggl[1-\sqrt{1-\frac{4R}{U^2} P_{\c,i} }\biggr],
\label{eq:P_b}
\end{align}
for $i=0,\ldots,N-1$. 
In the parallel configuration, the electric motor input-output map is given by
\begin{equation}
\label{eq:h_k}
\begin{split}
P_{\c,i} & = h_i(P_{\el, i},\omega_{\el, i}), \\
&= \kappa_{2}(\omega_{\el, i})P_{\el, i}^2 + \kappa_{1}(\omega_{\el, i})P_{\el, i} + \kappa_{0}(\omega_{\el, i}), 
\end{split}
\end{equation}
while for the series configuration we have
\begin{equation}
\label{eq:h_k2}
\begin{split}
P_{\e,i} & = h_i(P_{\drv, i},\omega_{\drv, i}),  \\
& = \kappa_{2}(\omega_{\drv, i})P_{\drv, i}^2  + \kappa_{1}(\omega_{\drv, i})P_{\drv, i} + \kappa_{0}(\omega_{\drv, i}),
\end{split}
\end{equation}
and
\begin{equation}
\label{eq:h_k3}
\begin{split}
P_{\gt,i} & = h_{\gen,i}(P_{\gen, i},\omega_{\gen, i}), \\
& = \nu_{2}(\omega_{\gen, i})P_{\gen, i}^2  + \nu_{1}(\omega_{\gen, i})P_{\gen, i} + \nu_{0}(\omega_{\gen, i}). 
\end{split}
\end{equation}
The aircraft dynamics are given in discrete time by
\begin{gather}
m_{i} v_{i} \Delta_{i} \gamma + m_{i} g \cos (\gamma_{i}) =   T_i \sin (\alpha_i) + \tfrac{1}{2}C_L(\alpha_{i})\rho S v_{i}^2
\label{eq:vertical_k}
\\
P_{\drv,i} = \tfrac{1}{2} m_{i} \Delta_{i} (v^2)  +  m_{i} g \sin (\gamma_{i}) v_{i}+\tfrac{1}{2}C_D(\alpha_{i})\rho S v_{i}^3
\label{eq:drive}
\end{gather}
for $i=0,\ldots,N-1$, where 
\[
\Delta_{i} (v^2) = (v^2_{i+1}   -   v^2_{i})/\delta , \quad
\Delta_{i} \gamma = (\gamma_{i+1}   -   \gamma_{i})/\delta.
\]
%
The power balance in discrete time for the parallel and series case respectively is given by
\begin{align}
\label{eq:P_parallel}
P_{\drv,i} &= P_{\gt,i} + P_{\el,i} , 
\\
\label{eq:P_series}
P_{\e,i} &= P_{\c,i} + P_{\gen,i} . 
\end{align}

\subsection{Parallel architecture}

For the parallel architecture the problem solved at the $k$th time step is
\begin{align}
& \min_{\substack{P_{\gt},\,P_{\el},\,P_{\drv},\,m,\\E,\,\omega_{\gt},\,\omega_{\el},\, \alpha}}
\quad \sum^{N-1}_{i=0} f_i(P_{\gt,i}, \omega_{\gt,i})\delta
\label{eq:min} \\
& \text{ s.t.}\quad 
\begin{aligned}[t]
& P_{\drv,i} = P_{\gt,i} + P_{\el,i} 
\\
&   P_{\drv,i} = \tfrac{1}{2} m_{i} \Delta_{i} v^2   +  m_{i} g \sin{(\gamma_{i})}v_{i} +\tfrac{1}{2}C_D(\alpha_{i})\rho S v_{i}^3 
\\
&  m_{i} v_{i} \Delta_{i} \gamma + m_{i} g \cos{\gamma_{i}} =   T_i \sin{\alpha_i} + \tfrac{1}{2}C_L(\alpha_{i})\rho S v_{i}^2 
\\
&   m_{i+1}  =m_{i} - f_i(P_{\gt,i}, \omega_{\gt,i}) \,\delta
\\
&   E_{i+1} = E_{i} - g_i \left(h_i\left(P_{\el, i},\omega_{\el, i}\right)\right)  \, \delta
\\
&   m_{0} = m(k\delta) 
\\
&   E_{0} = E(k\delta) 
\\
&  \underline{E} \leq E_{i} \leq  \overline{E}
\\
&  \underline{P}_{\gt} \leq P_{\gt,i} \leq  \overline{P}_{\gt}
\\
&  \underline{\omega}_{\gt} \leq \omega_{\gt,i} \leq  \overline{\omega}_{\gt}
\\
&  {\underline{P}_{\el}} \leq P_{\el,i} \leq  \overline{P}_{\el}
\\
&  {\underline{\omega}_{\el}} \leq \omega_{\el,i} \leq  \overline{\omega}_{\el}
\end{aligned}\nonumber
\end{align}
%
\subsection{Series architecture}

For the series architecture, the problem solved at the $k$th time step is
\begin{equation}
\min_{\substack{P_{\gt},\,P_{\e},\,P_{\drv},\,P_{\gen},\,P_{\c},\\m,\,E,\,\omega_{\gt},\,\omega_{\drv},\, \alpha}}
\quad \sum^{N-1}_{i=0} f_i(P_{\gt,i}, \omega_{\gt,i})\delta
\label{eq:min_series} 
\end{equation}
\vspace{-\baselineskip}
\begin{align*}
\quad \text{ s.t.} \quad
& P_{\e,i} = P_{\c,i} + P_{\gen,i} 
\\
& P_{\drv,i} = \tfrac{1}{2} m_{i} \Delta_{i} v^2   +  m_{i} g \sin{(\gamma_{i})}v_{i} +\tfrac{1}{2}C_D(\alpha_{i})\rho S v_{i}^3 
\\
& m_{i} v_{i} \Delta_{i} \gamma + m_{i} g \cos{\gamma_{i}} =   T_i \sin{\alpha_i} + \tfrac{1}{2}C_L(\alpha_{i})\rho S v_{i}^2
\\
&   m_{i+1}  =m_{i} - f_i(P_{\gt,i}, \omega_{\gt,i}) \,\delta
\\
&   E_{i+1} = E_{i} -  g_i (P_{\c,i})  \, \delta
\\
& P_{\e,i}  = h_i(P_{\drv, i},\omega_{\drv, i}) 
\\
& P_{\gt,i} = h_{\gen,i}(P_{\gen, i},\omega_{\gen, i})
\\
&   m_{0} = m(k\delta) 
\\
&   E_{0} = E(k\delta) 
\\
&  \underline{E} \leq E_{i} \leq  \overline{E}
\\
&  \underline{P}_{\gt} \leq P_{\gt,i} \leq  \overline{P}_{\gt}
\\
&  \underline{\omega}_{\gt} \leq \omega_{\gt,i} \leq  \overline{\omega}_{\gt}
\\
&  {\underline{P}_{\drv}} \leq P_{\drv,i} \leq  \overline{P}_{\drv}
\\
&  {\underline{\omega}_{\drv}} \leq \omega_{\drv,i} \leq  \overline{\omega}_{\drv}
\end{align*}
%



\section{Convex relaxation}\label{sec:convex}

The optimisation problems in (\ref{eq:min}) and (\ref{eq:min_series}) are nonconvex, which makes a real-time implementation of an MPC algorithm that relies on its solution computationally intractable. In this section a convex approximation is proposed that is suitable for an online solution. We make three simplifications: 1) we prescribe a flight profile and impose an assumption on the monotonicity of the loss map functions which results in convex loss map functions and allows their coefficients to be computed a priori; 2) we reformulate the dynamics as a quadratic function of aircraft mass under mild assumptions; 3) we introduce a \mdsedittris{lossless} change of  optimisation variables that shifts the nonlinear term in the battery update equation to the power balance inequality.

\subsection{Reformulation of the loss map functions}

We assume that the aircraft speed $v_i$ and flight path angle $\gamma_i$ are chosen externally by a suitable guidance algorithm for $i=0,\ldots,N-1$. \mdseditbis{This assumption is reasonable for an actual air traffic management application where flight corridors are prescribed.}
For the series configuration, we assume 
\mdseditbis{that the generator speed is constant:  $\omega_{\gen,i} = \omega_{\gen}^{\ast}$, $\forall i$, where the optimal speed $\omega_{\gen}^{\ast}$ is determined empirically so as to operate the turbo generator set at its maximum efficiency. This allows us to fix the coefficients in \eqref{eq:h_k3} to constant values and express $h_{\gen,i}(P_{\gen,i},\omega_{\gen,i})$ as a convex quadratic function of $P_{\gen,i}$
\begin{equation}
h_{\gen,i}(P_{\gen,i}) = \nu_{2}P_{\gen,i}^2  + \nu_{1}P_{\gen,i}  + \nu_{0},
\label{eq:h_P3}
\end{equation}
with constant coefficients $\nu_{2} \geq 0$, $\nu_{1} >0$. }

\mdseditbis{For the parallel configuration, we assume for simplicity that the gas turbine, electric motor and fan share a common shaft rotation speed, i.e.\ $\omega_{\gt,i} = \omega_{\el,i} = \omega_{\drv,i}$, $\forall i$.   If the fan shaft speed is known at each time step of the prediction horizon, then the coefficients in
\eqref{eq:h_k} and \eqref{eq:h_k2} can be estimated from a set of polynomial approximations of $h_i(\cdot)$ at a pre-determined set of speeds. 
This allows $h_i(P_{\el,i},\omega_{\el,i})$ and $h_i(P_{\drv,i},\omega_{\drv,i})$
to be replaced by time-varying convex functions of power alone
\begin{align}
h_i(P_{\el,i}) = \kappa_{2,i}P_{\el,i}^2  + \kappa_{1,i}P_{\el,i}  + \kappa_{0, i},
\label{eq:h_P}
\\
 h_i(P_{\drv,i}) = \kappa_{2,i}P_{\drv,i}^2  + \kappa_{1,i}P_{\drv,i}  + \kappa_{0, i},
\label{eq:h_P2}
\end{align}
with $\kappa_{2,i} \geq 0$, $\kappa_{1,i} >0$, $\kappa_{2,i} \geq 0$, $\kappa_{1,i} >0$,
for all $i$. Regarding the gas turbine fuel map, since the spool speed is assumed constant, the coefficients are independent of gas turbine spool speed such that $f_i(P_{\gt,i},\omega_{\gt,i})$ can also be replaced by a convex functions of power alone
\begin{equation}
f_i(P_{\gt,i}) = \beta_{2} P_{\gt,i}^2 + \beta_{1} P_{\gt,i}  + \beta_{0}, 
\label{eq:f_P}
\end{equation}
with $\beta_{2} \geq 0$, $\beta_{1} > 0$.}

If moreover we assume that these functions are non-decreasing \mdseditbis{as suggested in \cite{east2018}}, the following hold: $P_{\el,i} \geq -\kappa_{1,i}/2\kappa_{2,i}$, $P_{\drv,i} \geq -\kappa_{1,i}/2\kappa_{2,i}$, $P_{\gen,i} \geq -\nu_{1}/2\nu_{2}$, $P_{\gt,i} \geq -\beta_{1}/2\beta_{2}$ for all $i$, which requires  new lower bounds on power. In the parallel configuration, the new bounds are given by
\begin{align}
\label{eq:pow1}
\underline{P}_{\el,i} = \max{ \left\{ \underline{P}_{\el}, - \frac{\kappa_{1,i}}{2\kappa_{2,i}} \right\} }, 
\\
\label{eq:pow2}
\underline{P}_{\gt} = \max{ \left\{ \underline{P}_{\gt}, - \frac{\beta_{1}}{2\beta_{2}} \right\} }, 
\end{align}
whereas in the series configuration only the gas turbine bound should be updated, as follows
\begin{equation}
\label{eq:pow3}
\underline{P}_{\gt} = \max{ \left\{ \underline{P}_{\gt}, - \frac{\beta_{1}}{2\beta_{2}} , h_{\gen,i}\left( - \frac{\nu_{1}}{2\nu_{2}}\right) \right\} }, 
\end{equation}
since we can enforce the monotonicity condition on the drive power a priori when prescribing the drive power profile. 

In order to estimate the shaft speed $\omega_{\drv,i}$, and hence determine the coefficients in (\ref{eq:h_P})-(\ref{eq:f_P}), we use a pre-computed look-up table relating the drive power to rotational speed of the fan, for a given altitude,  Mach number, and air conditions (temperature and specific heat at constant pressure). This enables the shaft speed to be determined as a function of the fan power output at each discrete-time step along the flight path. Although $P_{\drv,i}$ depends on the aircraft mass $m_i$, which is itself an optimisation variable, a prior estimate of the required power output can be obtained 
by assuming a constant mass $m_i=m_0$ for all $i$. It was shown in \cite{doff2020optimal} that this assumption has a negligible effect on solution accuracy.

Note that since the rotation speeds are prescribed, all constraints on shaft rotation speeds can be removed from the optimisation (and checked a priori).  The same remark holds for the constraints on the drive power. 

\subsection{Reformulation of the dynamics}

To express the dynamics in a form suitable for convex programming, we simplify the dynamical equations and combine the equations that constrain the aircraft motion  as follows.
First we express the drag and lift coefficients, $C_D$ and $C_L$, as functions of the angle of attack $\alpha$. Over a restricted domain and for given Reynolds and Mach numbers, the drag and lift coefficients can be expressed respectively as a quadratic non-decreasing function and a linear non-decreasing function \cite{abbott1945}
\begin{alignat}{2}
C_D(\alpha_{i}) &= a_2 \alpha_{i}^2 + a_1 \alpha_{i} + a_0, &\qquad  a_2 &> 0, 
\label{eq:C_D}
\\
C_L(\alpha_{i}) &= b_1 \alpha_{i} + b_0, &\qquad  b_1 &> 0,
\label{eq:C_L}
\end{alignat}
for $\underline{\alpha} \leq \alpha_i \leq \overline{\alpha}$.
Secondly, assuming that the contribution of the thrust in the vertical direction is negligible%
\footnote{\mdsedittris{This assumption was checked in simulations, where it was found that the solution satisfies $\alpha < 2^\circ$, which supports this assumption.}}, the term $T \sin{(\alpha)}$ can be neglected from \eqref{eq:vertical_k}.
Finally, combining \eqref{eq:vertical_k}, \eqref{eq:drive}, \eqref{eq:C_D} and \eqref{eq:C_L}, the angle of attack can be eliminated from the expression for $P_{\drv,i}$, which can be expressed as a quadratic function of the aircraft mass, $m_{i}$, as follows
\begin{equation}\label{eq:P_drv}
    P_{\drv,i} = \eta_{2,i} m_{i}^2 + \eta_{1,i} m_{i} + \eta_{0,i} , 
\end{equation}
where
\begin{align*}
&    \eta_{2,i} = \frac{ 2 a_2 (v_{i} \Delta_{i} \gamma + g \cos{\gamma_{i}}  )^2 }{ b_1^2 \rho S v_{i}} ,
\\
&    \eta_{1,i} = \tfrac{1}{2}\Delta_{i} v^2 + g \sin{\gamma_{i}} v_{i}   - \frac{2 a_2 b_0  (v_{i} \Delta_{i} \gamma
    + g \cos{\gamma_{i}} ) v_{i}}{b_1^2} 
\\
&\qquad\quad    + \frac{a_1}{ b_1 }  ( v_{i} \Delta_{i} \gamma + g \cos{\gamma_{i}}  ) v_{i} ,
\\
&    \eta_{0,i} = \tfrac{1}{2} \rho S v_{i}^3 \Bigl(  \frac{a_2 b_0^2}{b_1^2} - \frac{a_1 b_0}{b_1} + a_0  \Bigr). 
\end{align*}
Since the flight path angle $\gamma_i$ and speed $v_i$ are determined a priori, the coefficients $\eta_{0,i}$, $\eta_{1,i}$, $\eta_{2,i}$ are fixed. Moreover $\eta_{2,i} > 0$ for all $i$, so the drive power is a convex function of  $m_{i}$. Note that there is no guarantee that satisfaction of \eqref{eq:P_drv} enforces \eqref{eq:vertical_k} and \eqref{eq:drive} individually.  In practice,  assuming that we have full control over the eliminated variable $\alpha$ (via the elevator and fans), both individual dynamical equations can be satisfied a posteriori. \finaledit{The existence of a faster inner flight control loop on angle of attack and thrust ensures that the correct trajectory is followed}. The bounds  $\underline{\alpha}\leq\alpha_i\leq\overline{\alpha}$ need to be checked a posteriori.

\subsection{Reformulation of power balance}
Let the rate of change of fuel mass be $\varphi_i := f_i(P_{\gt,i})$. Using this new variable, the power balance can be enforced by
\[
\varphi_i = f_{\varphi, i} (m_i, P_{\b,i}), 
\]
where the function $f_{\varphi,i}$ is defined  
\begin{equation*}
  f_{\varphi,i}=\begin{cases}
    f_i\left(P_{\drv,i}(m_i) - h_i^{-1}\left(g_i^{-1}\left(P_{\b,i}\right)\right)  \right) & \text{if $\mathcal{T} = \mathcal{P}$}\\
    f_i\left( h_{\gen,i} \left( h_i\left(P_{\drv,i}(m_i)\right) - g_i^{-1}\left(P_{\b,i}\right) \right) \right)  & \text{if $\mathcal{T} = \mathcal{S}$}
    \end{cases}
\end{equation*}
where $P_{\drv,i}$ is given by equation \eqref{eq:P_drv}. This formulation unifies the treatment of series and parallel configurations and eliminates the variables $P_{\e}$, $P_{\gen}$ and $P_{\el}$ from the optimisation problem. Moreover, since the functions $f_i(\cdot)$, $h_i(\cdot)$, $h_{\gen,i}(\cdot)$, $g_i(\cdot)$, $P_{\drv,i}(\cdot)$ are convex, twice differentiable, non-decreasing, one-to-one functions,  the function $f_{\varphi,i}(\cdot)$ is also convex.

We can construct a convex program by relaxing the power balance equality to the inequality
\begin{equation}
  \varphi_i \geq f_{\varphi, i} (m_i, P_{\b,i}) ,
\end{equation}
\finaledit{which is necessarily satisfied with equality at the optimum since the form of the objective in \eqref{eq:min} and \eqref{eq:min_series} ensures that any feasible solution that does not satisfy this constraint with equality is suboptimal.}

The constraints on gas turbine power and electric motor power are replaced by constraints on rate of change of fuel mass and on battery power, respectively,
\begin{gather}
\underline{\varphi}_i  \leq \varphi_i \leq  \overline{\varphi}_i,
\\
{\underline{P}_{\b,i}} \leq P_{\b,i} \leq  \overline{P}_{\b,i}, 
\end{gather}
with $\underline{\varphi}_i = f_i(\underline{P}_{\gt})$, $\overline{\varphi}_i = f_i(\overline{P}_{\gt})$. 
Here 
\[
\underline{P}_{\b,i}  = \begin{cases}
    g_i(h_i(\underline{P}_{\el,i})) & \text{if $\mathcal{T} = \mathcal{P}$} \\
   \widetilde{f}_i(\underline{P}_\el, \overline{P}_{\gt} )  & \text{if } \mathcal{T} = \mathcal{S}\end{cases}
\]
where $\widetilde{f}_i(x,y):=g_i( h_i(x)-h_{\gen,i}^{-1} (y) )$.
Furthermore, to ensure that  $g_i(\cdot)$ is real-valued we require $P_{\c,i} \leq U^2/4R$, and hence 
\[   
\overline{P}_{\b,i}  = \begin{cases}
    g_i(h_i(\overline{P}_{\el,i})) & \text{if } \mathcal{T} = \mathcal{P} \\
     \min \bigl\{ \widetilde{f}_i( \overline{P}_\el, \underline{P}_{\gt}),\tfrac{U^2}{2R} \bigr\}  & \text{if } \mathcal{T} = \mathcal{S}
     \end{cases}
\]
where $\overline{P}_{\el,i}$ is defined for the parallel configuration by
\[
\overline{P}_{\el,i} = \min{\{ \overline{P}_{\el}, r_{\text{max},i} \}}  
\]
with $r_{\text{max},i}=\max {\{ x: 1 - 4R/U^2 h_{i}(x)= 0\}}$.
\subsection{Convex program}
A unified convex program can thus be formulated as follows 
\begin{align}
& \min_{\substack{\varphi,\,P_{\b},\,m,\,E}}
\quad \sum^{N-1}_{i=0} \varphi_i\delta
\label{eq:min2} \\
\quad & \text{ s.t.} \quad 
\begin{aligned}[t]
& \varphi_i \geq f_{\varphi,i}\left(m_{i}, P_{\b,i}\right)
\\
&   m_{i}  =m(k\delta) - \sum_{l=0}^{i-1} \varphi_l \,\delta
\\
&   E_{i} = E(k\delta) -  \sum_{l=0}^{i-1} P_{\b,l}  \, \delta
\\
&  \underline{E} \leq E_{i} \leq  \overline{E}
\\
&  \underline{\varphi}_i  \leq \varphi_i \leq  \overline{\varphi}_i
\\
&  {\underline{P}_{\b,i}} \leq P_{\b,i} \leq  \overline{P}_{\b,i}
\end{aligned}
\nonumber
\end{align}
where the bounds $\underline{\varphi}_i$, $\overline{\varphi}_i$, $\underline{P}_{\b,i}$, $\overline{P}_{\b,i}$ are given by
\[
\overline{\varphi}_i = f_i(\overline{P}_{\gt}),
\]
and, for $\mathcal{T} = \mathcal{P}$:
\begin{align*}
\underline{\varphi}_i &= 
\max\bigl\{ f_i( \underline{P}_{\gt} ),f_i(  - \tfrac{\beta_{1}}{2\beta_{2}} ) \bigr\}
\\
\underline{P}_{\b,i}  &= 
\max\bigl\{ g_i\bigl(h_i(\underline{P}_{\el})\bigr), g_i\bigl(h_i(- \tfrac{\kappa_{1,i}}{2\kappa_{2,i}})\bigr) \bigr\}
\\
\overline{P}_{\b,i}  &=
\min\bigl\{  g_i\bigl(h_i(\overline{P}_{\el})\bigr),  g_i\bigl(h_i(r_{\text{max},i})\bigr) \bigr\},
\end{align*}
and, for $\mathcal{T} = \mathcal{S}$:
\begin{align*}
\underline{\varphi}_i &=
\max\bigl\{ f_i( \underline{P}_{\gt} ),f_i(  - \tfrac{\beta_{1}}{2\beta_{2}} ), f_i\bigl( h_{\gen,i}( - \tfrac{\nu_{1}}{2\nu_{2}}) \bigr) \bigr\} 
\\
\underline{P}_{\b,i}  &= 
    \widetilde{f}_i ( \underline{P}_\el, \overline{P}_{\gt})
\\
\overline{P}_{\b,i}  &=
\begin{aligned}[t] \min \Bigl\{ 
\widetilde{f}_i(\overline{P}_\el,\underline{P}_{\gt}),
&\widetilde{f}_i(\overline{P}_\el,- \tfrac{\beta_{1}}{2\beta_{2}} ),
\\
&\widetilde{f}_i\bigl( \overline{P}_\el, h_{\gen,i}( - \tfrac{\nu_{1}}{2\nu_{2}}) \bigr),  \tfrac{U^2}{2R} \Bigr\}.
\end{aligned}
\end{align*}
\section{Alternating Direction Method of Multipliers} \label{sec:admm}
If $\underline{E} \leq E_{i}  - \overline{P}_{\b,i}\, \delta \leq  \overline{E}$ $\forall i$, \finaledit{so that at each time step there is enough energy in the battery to operate the electric motor at its maximum capacity}, then the solution of  \eqref{eq:min2} is given trivially by $P_{\b,i}^{\ast} = \overline{P}_{\b,i} $, $\forall i$, for both architectures. If this condition is not satisfied, then an optimisation scheme is needed to solve problem \eqref{eq:min2}.
To make real-time implementation possible we propose a specialised ADMM algorithm \cite{boyd11}.
Problem \eqref{eq:min2} can be equivalently stated with inequality constraints appended to the objective function using indicator functions $\Lambda^x(x)$,
\begin{alignat}{3}
& \min_{\substack{\varphi,\,P_{\b},\,m,\\E,\, \chi,\, \xi,\, \zeta}} & & \sum^{N-1}_{i=0}  \xi_i \delta + \Lambda^{\chi}(\chi_i) + \Lambda^E(E_i) + \Lambda^{\varphi}(\varphi_i) + \Lambda^{P_{\b}}(P_{\b,i})
\label{eq:min_dummy} \\
& \ \text{ s.t.} \quad & &  \chi_i = \xi_i -  f_{\varphi,i} (m_{i}, P_{\b,i})
\nonumber \\
& & &  m_{i}  =m(k\delta) - \sum_{l=0}^{i-1} \xi_l \delta 
\nonumber \\
& & &   E_{i} = E(k\delta) - \sum_{l=0}^{i-1}\zeta_l \delta 
\nonumber \\
& & &   \xi_i = \varphi_i
\nonumber \\
& & &   \zeta_i = P_{\b,i} 
\nonumber 
\end{alignat}
with $\underline{\chi} = 0$, $\overline{\chi} = \infty$, and 
\[
\Lambda^x(x)= \begin{cases}
    0 & \text{if $\underline{x} \leq x \leq \overline{x}$},\\
   \infty & \text{otherwise}.
  \end{cases}
\]
Note that we have introduced dummy variables $\xi$ and $\zeta$ in order to simplify the solver iterations by separating variables.

We define an augmented Lagrangian function as
\[
\label{eq:lagrange}
\begin{split}
 & L(\chi,\xi,\zeta,E,P_{\b},\varphi,m,\lambda_1,\lambda_2,\lambda_3,\lambda_4,\lambda_5) = \\ 
 &\sum_{i=0}^{N-1}\bigl( \xi_i \delta + \Lambda^{\chi}(\chi_i) + \Lambda^{\!E}(E_i) + \Lambda^{\varphi}(\varphi_i) + \Lambda^{\!P_{\b}}(P_{\b,i}) \bigr) \\
& + \frac{\sigma_1}{2} \sum_{i=0}^{N-1} \bigl(\chi_i - \xi_i + f_{\varphi,i} (m_{i}, P_{\b,i}) + \lambda_{1,i} \bigr)^2  \\
& + \frac{\sigma_2}{2} \| m - m(k\delta)\Phi  + \Psi \xi +\lambda_2 \|^2 \\
& + \frac{\sigma_3}{2} \| E - E(k\delta)\Phi + \Psi \zeta + \lambda_3 \|^2 \\
& + \frac{\sigma_4}{2} \| \xi -\varphi + \lambda_4 \|^2 \\
& + \frac{\sigma_5}{2} \| \zeta - P_{\b} + \lambda_5 \|^2,
\end{split}
\]
where $\lambda_i$ is a Lagrange multiplier and $\sigma_i$ is a penalty parameter associated with the $i$th constraint, $\Phi$ is a vector of ones, and $\Psi$ is the strictly lower triangular matrix with zeros on the diagonal and all other lower triangular elements equal to~$\delta$. 

Problem \eqref{eq:min_dummy} can be rearranged in the canonical form
\begin{alignat}{2}
& \min_{x,\,z}
& \quad & \hat{f}(x) + \hat{g}(z)
\\
& \text{ s.t.} & & b(z) + Bx = c 
\nonumber
\end{alignat}
with 
\begin{gather*}
x=\begin{bmatrix} \chi^\top & \xi^\top & \zeta^\top & E^\top & \varphi^\top \end{bmatrix}^\top, \quad 
z=\begin{bmatrix} m^\top & P_b^\top \end{bmatrix}^\top, 
\\
\lambda = \begin{bmatrix} \lambda_1^\top & \lambda_2^\top & \lambda_3^\top & \lambda_4^\top & \lambda_5^\top \end{bmatrix}^\top, 
\\
\hat{f}(x) = \sum^{N-1}_{i=0}  \xi_i \delta + \Lambda^{\chi}(\chi_i) + \Lambda^E(E_i) + \Lambda^{\varphi}(\varphi_i) , 
\\
\hat{g}(z) =  \sum^{N-1}_{i=0}  \Lambda^{P_{\b}}(P_{\b,i}), 
\\
B=\begin{bmatrix} I & -I & 0 & 0 & 0\\
0 & \Psi & 0 & 0 & 0\\
0 & 0 & I & 0 & 0 \\
0 & 0 & \Psi & I & 0\\ 
0 & I & 0 & 0 & -I\end{bmatrix}, \quad 
b(z) = \begin{bmatrix} f_\varphi (m, P_b) \\ m \\-P_{b}\\ 0 \\  0 \end{bmatrix}, 
\\
c=\begin{bmatrix} 0 & \Phi^\top m(k\delta)& 0 & \Phi^\top E(k\delta)& 0 \end{bmatrix}^\top.
\end{gather*}
We define the primal and dual residuals $r^{j+1} = b(z^{j+1}) + Bx^{j+1} - c$ and $s^{j+1} = [\nabla_{z}b (z^{j+1})]^\top R^j B (x^j - x^{j+1}) $, where $R^j= \text{diag}(\sigma_1^j I, \sigma_2^j I, \sigma_3^j I, \sigma_4^j I, \sigma_5^j I)$. Note that $0$ and $I$ are compatible zero and identity matrices. \mdsedit{By comparison with  \cite{boyd11}, the present algorithm deals with a nonlinear $b(z)$ function in the equality constraint, which requires that the dual residual is  defined in terms of the Jacobian $\nabla_{z}b$.} 

The ADMM iteration update is given by
\begin{align*}
&\chi^{j+1} = \proj^{\chi} ( \xi^j - f^j_\varphi - \lambda_{1}^j ),
\\
&\xi^{j+1} = \bigl( (\sigma^j_1\!+\! \sigma^j_4)I + \sigma^j_2 \Psi^\top\! \Psi \bigr)^{-1}\! \Bigl[ -\Phi\delta + \sigma^j_1(\chi^{j+1} + f^j_\varphi + \lambda_1^j)
\\
&\qquad\qquad - \sigma^j_2 \Psi^\top (m^j - m(k\delta)\Phi + \lambda_2^j) + \sigma^j_4 \bigl( \varphi^j - \lambda_4^j\bigr) \Bigr], 
\\
& 
\zeta^{j+1} = (\sigma^j_5 I + \sigma^j_3 \Psi^\top\Psi)^{-1} \Bigl[
\begin{aligned}[t] 
& -\sigma^j_3\Psi^\top\bigl( E^j - E(k\delta)\Phi + \lambda_3^j \bigr) 
\\
& + \sigma^j_5 (P_{\b}^j - \lambda_5^j) \Bigr], 
\end{aligned}
\\
& E^{j+1} = \proj^{E} \bigl( E(k\delta)\Phi - \Psi\zeta^{j+1} - \lambda_{3}^j \bigr), 
\\
& P_{\b,i}^{j+1} = \proj^{P_{\b}} \Bigl( \begin{aligned}[t]
& \arg\min_{P_{\b,i}} \Bigl\{
\frac{\sigma^j_1}{2} \bigl[ \chi_i^{j+1} - \xi_i^{j+1} + f_{\varphi,i}(m_{i}^j, P_{\b,i})\\
& + \lambda_{1,i}^j \bigr]^2 + \frac{\sigma^j_5}{2} \bigl[ \zeta_i^{j+1} - P_{\b,i} + \lambda_{5,i}^j \bigr]^2
\Bigr\}\Bigr), 
\end{aligned}
\\
& \varphi^{j+1} = \proj^{\varphi} \bigl( \xi^{j+1} + \lambda_4^j  \bigr) , 
\\
& m_i^{j+1} = \begin{aligned}[t] 
&\arg\min_{m_i} \Bigl\{ \frac{\sigma^j_1}{2} \Bigl[ \chi_i^{j+1} \!- \xi_i^{j+1} \!+ f_{\varphi,i} (m_{i}, P_{\b,i}^{j+1}) + \lambda_{1,i}^j \Bigr]^2 
\\
& +  \frac{\sigma^j_{2}}{2} \Bigl[ m_i - m(k\delta)\Phi + [\Psi \xi^{j+1}]_i + \lambda_2^j\Bigr]^2 \Bigr\}, 
\end{aligned}
\\
& \lambda_1^{j+1} = \lambda_1^j + \chi^{j+1} - \xi^{j+1} + f_\varphi^{j+1}, 
\\
& \lambda_2^{j+1} = \lambda_2^j + m^{j+1} - m(k\delta)\Phi + \Psi \xi^{j+1}, 
\\
& \lambda_3^{j+1} = \lambda_3^j + E^{j+1} - E(k\delta)\Phi + \Psi \zeta^{j+1}, 
\\
& \lambda_4^{j+1} = \lambda_4^j + \xi^{j+1} - \varphi^{j+1}, 
\\
& \lambda_5^{j+1} = \lambda_5^j + \zeta^{j+1} - P_{\b}^{j+1}, 
\end{align*}
where $f^j_\varphi=\smash{[ f_{\varphi,0}^j \ \cdots \ f_{\varphi,{N-1}}^j]^\top}$ and $\pi^x (y)$ denotes the projection  $\max\{\min\{ y, \overline{x}\}, \underline{x}\}$. The penalty parameters $\sigma_n^j$, $n=1,2,3,4,5$ are updated at intervals of $F_\sigma$ iterations (provided $10<\max{\left\{\frac{||r^{j+1}||}{\max{\{ || b(z^{j+1})||, ||B x^{j+1}||, ||c||\}}},\frac{||s^{j+1}||}{|| \nabla_{z}b (z^{j+1})^\top \lambda^{j+1} ||} \right\}} $) according to the rule
\begin{align*}
\tau^{j+1} &=  \begin{cases}
    \Gamma & \text{if $ 1 \leq \Gamma < \tau_{\text{max}} $},\\
    \Gamma^{-1} & \text{if $ \tau_{\text{max}}^{-1} < \Gamma < 1 $}, \\
    \tau_{\text{max}}  & \text{otherwise},  
  \end{cases}
\\
\sigma^{j+1}_n &=  \begin{cases}
    \sigma^{j}_n \tau^{j+1} & \text{if $ \left\lVert r^{j+1}_n \right \rVert > \mu \left\lVert s^{j+1} \right \rVert $},\\
    \sigma^{j}_n/\tau^{j+1} & \text{if $ \left\lVert s^{j+1} \right \rVert > \mu \left\lVert r^{j+1}_n \right \rVert $}, \\
    \sigma^{j}_n  & \text{otherwise},  
  \end{cases}
\\
R^{j+1} &= \diag(\sigma_1^{j+1}I, \sigma_2^{j+1}I, \sigma_3^{j+1}I, \sigma_4^{j+1}I, \sigma_5^{j+1}I), 
\end{align*}
where $\Gamma = \sqrt{\left\lVert r^{j+1} \right \rVert / \left\lVert s^{j+1} \right \rVert}$   and $r_n^{j+1}$ denotes the rows of $r^{j+1}$ associated with the $n$th constraint, $1\leq n \leq 5$.

The updates for $\chi$, $E$ and the multipliers $\lambda_1$, $\lambda_2$, $\lambda_3$, $\lambda_4$,  $\lambda_5$ involve only vector additions, summations and projections. 
%
%
%
\finaledit{The equations $\bigl((\sigma^j_1 + \sigma^j_4) I + \sigma^j_3 \Psi^\top\Psi\bigr)\xi = c_1$ and $(\sigma^j_5 I + \sigma^j_3 \Psi^\top\Psi)\zeta = c_2$ can be solved for $\xi$ and $\zeta$ (for given $c_1$ and $c_2$) in $O(N)$ operations using appropriate Cholesky factorisations (for details see~\cite{east2020ev}, Prop.~3). The Cholesky factors can be reused until the penalty parameters $\sigma^j_n$ are updated, otherwise the updates for $\xi$ and $\zeta$ require only scalar multiplications and vector summations.}
The updates for $P_b$ and $m$ require minimisation of scalar convex functions and can be performed using Newton's method. 


The algorithm is initialised with
 \[
 \begin{split}
 & P_{b}^{0}= \Phi \overline{P}_{b},  \quad \zeta^0 = P_{b}^{0}, \quad \xi^0= \underline{\varphi}, \quad \varphi^0= \xi^0, \\
 & E^0 = \pi^E \left (  \Phi E(k\delta)  - \mathrm{\Psi} \zeta^0 \right), \quad m^0  = \Phi m(k\delta) - \mathrm{\Psi} \xi^0 \\
 & \chi^{0} = \proj^{\chi} ( \xi^0 - f^0_\varphi), \quad \lambda_{1}^{0}=  \lambda_{2}^{0}= \lambda_{3}^{0}= \lambda_{4}^{0}= \lambda_{5}^{0}=0\Phi\\
 & R^0 = \text{diag}(50I, 3.69 \times 10^{-7} I,  6.96 \times 10^{-7}I, 20.29I, 0.83I), 
 \end{split}
\]
and stopped when  $\left\lVert r^{j+1} \right \rVert  \leq  \epsilon_P $ and $\left\lVert s^{j+1} \right \rVert \leq  \epsilon_D $ or $j> 10^{5} $, where, following \cite{boyd11},
\begin{align*}
\epsilon_P &= \sqrt{5N} \epsilon_{\text{abs}} + \epsilon_{\text{rel}} \max{\{ || b(z^{j+1})||, ||B x^{j+1}||, ||c||\}},
\\
\epsilon_D &= \sqrt{2N} \epsilon_{\text{abs}} + \epsilon_{\text{rel}} || [\nabla_{z}b (z^{j+1})]^\top \lambda^{j+1} ||.
\end{align*}
\mdseditbis{The penalty parameters $\sigma_n^j$ are  initialised so that all terms of the Lagrangian are initially of the same order of magnitude.}

\section{Numerical results}
\label{sec:simulation_results}

In this section we introduce an energy management case study involving a representative hybrid-electric passenger aircraft and solve optimisation problem \eqref{eq:min2} within this context using the ADMM algorithm as presented in section \ref{sec:admm}. The simulation results are analysed and the performance of the algorithm is discussed in terms of its computational requirements and robustness to variations in model parameters. 

\subsection{Simulation scenario}\label{sec:sim_setup}

The parameters of the model used in simulations are shown in Table \ref{tab:param}. These are based on published data for the BAe~146 aircraft. The conventional BAe 146 propulsion system is replaced by hybrid-electric propulsion systems\footnote{In order to maintain a constant MTOW, the excess mass from the batteries, electric motors, generators and electrical distribution systems can be compensated by cuts in passenger count and fuel mass.} in either parallel or series configuration (as illustrated in Figs.~\ref{fig:propulsion} and~\ref{fig:propulsion_series}), \mdsedittris{both of which were equipped with the same battery size}. The aircraft is powered by a combination of $4$ such systems.

\begin{table}[h!]
\centering
\begin{tabular}{llll}
\hline
\textbf{Parameter} & \textbf{Symbol} & \textbf{Value} & \textbf{Units} \\ \hline
Mass (MTOW)     &    $m$    &   $42000$    &    \si{kg}   \\ \hline
Gravity acceleration     &    $g$    &   $9.81$    &    \si{m.s^{-2}}   \\ \hline
Wing area &   $S$     &    $77.3$   &    \si{m^2}   \\ \hline
Density of air &   $\rho$     &    $1.225$   &    \si{kg.m^{-3}} \\ \hline
\multirow{2}*{Lift coefficients} &   $b_0$   &    $0.43$   &    \si{-}  \\
&   $b_1$     &    $0.11$   &    \si{deg^{-1}}  \\ \hline
 \multirow{3}*{Drag coefficients} &   $a_0$   &   $0.029$   &    \si{-}  \\
&   $a_1$   &   $0.004$   &    \si{deg^{-1}} \\
&   $a_2$   &   $5.3\mathrm{e}{-4}$   &    \si{deg^{-2}}\\ \hline
Angle of attack range &   $\left[\underline{\alpha}; \overline{\alpha}\right]$     &    $\left[-3.9; 10\right]$   &    \si{deg}   \\ \hline
$\#$ of propulsion systems &   $n$     &    $4$   &    \si{-} \\ \hline
Total fuel mass &   $m_{\fuel}$     &    $1000 \times n$   &    \si{kg}   \\ \hline
Total battery mass &   $m_{\b}$     &    $2000 \times n$   &    \si{kg}   \\ \hline
Battery energy density &   $e_{\b}$     &    $0.875$   &    \si{MJ.kg^{-1}}   \\ \hline
\multirow{2}*{Fuel map coefficients} &   $\beta_0$     &    $0.0327$    &   \si{kg.s^{-1}} \\ 
&   $\beta_1$     &    $0.0821$    &   $\si{kg.{MJ}^{-1}}$ \\  \hline
\multirow{2}*{Generator coefficients} &   $\nu_0$     &    $0.08$    &   \si{MW} \\ 
&   $\nu_1$     &    $1$    &   $\si{-}$ \\  \hline
Total battery SOC range &   $\left[\underline{E}; \overline{E}\right] \times n$     &    $\left[350; 1487\right] \times n$   &    \si{MJ} \\ \hline
Gas turbine power range &   $\left[\underline{P}_{\gt}; \overline{P}_{\gt}\right]$     &    $\left[0; 5\right]$  &    \si{MW}   \\ \hline
Motor power range &   $\left[\underline{P}_{\el}; \overline{P}_{\el}\right]$     &    $\left[0; 5\right]$   &    \si{MW}   \\ \hline
Battery o.c. voltage & $U$ & $1500$ & \si{V} \\ \hline
Battery resistance & $R$ & $0.035$ & \si{ohm} \\ \hline
Flight time &   $T$     &    $3600$   &    \si{s} \\ \hline
\end{tabular}
\caption{Model parameters.}
\label{tab:param}
\end{table}

\begin{figure}
\centerline{\includegraphics[width = .6\textwidth]{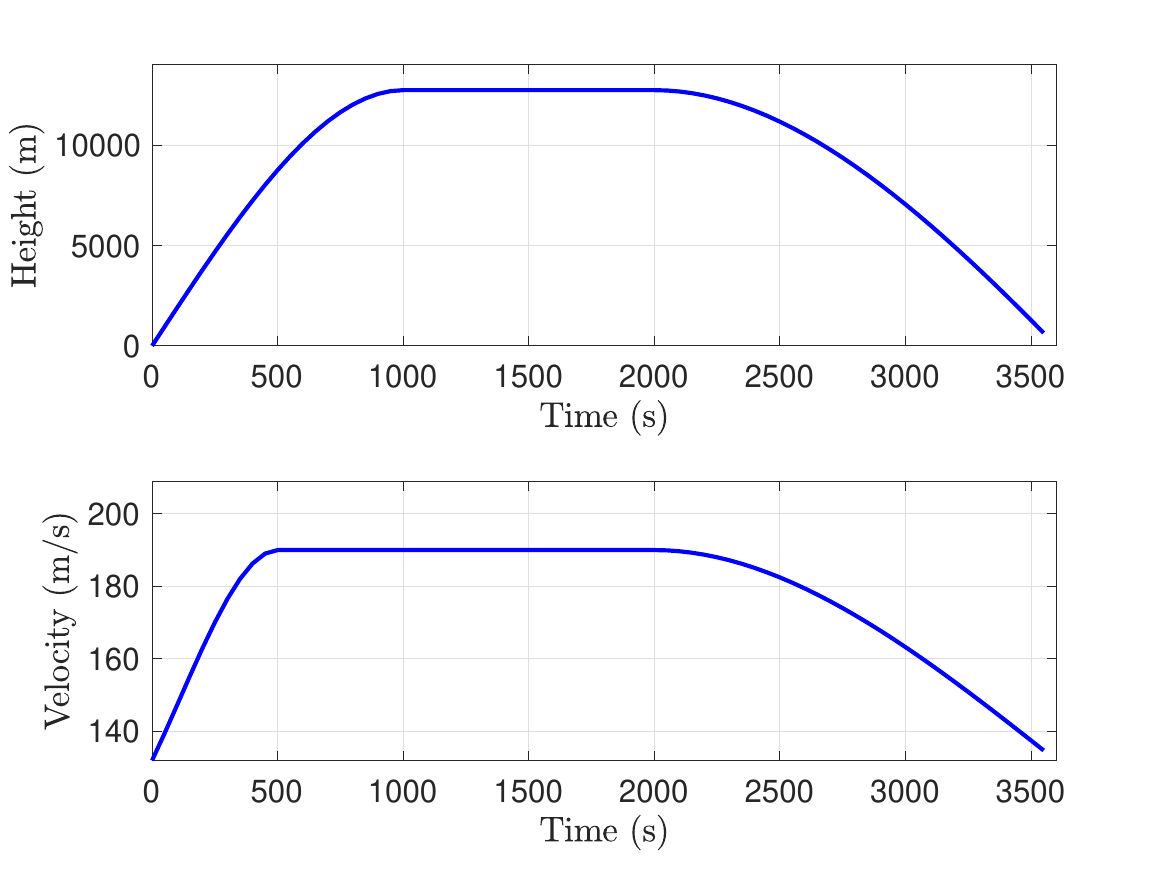}}
\caption{Height and velocity profiles for the mission.}
\label{fig:profile}
\end{figure}

For the purposes of this study it is assumed that velocity and height profiles are known a priori as a result of a fixed flight plan determined prior to take-off.  
We consider an exemplary 1-hour flight at a true airspeed (TAS) of $190$\,$\si{m/s}$ for a typical 100-seat passenger aircraft. The flight path (height and velocity profile) is shown in Figure \ref{fig:profile}. 

\begin{figure}
\centerline{\includegraphics[width = .6\textwidth]{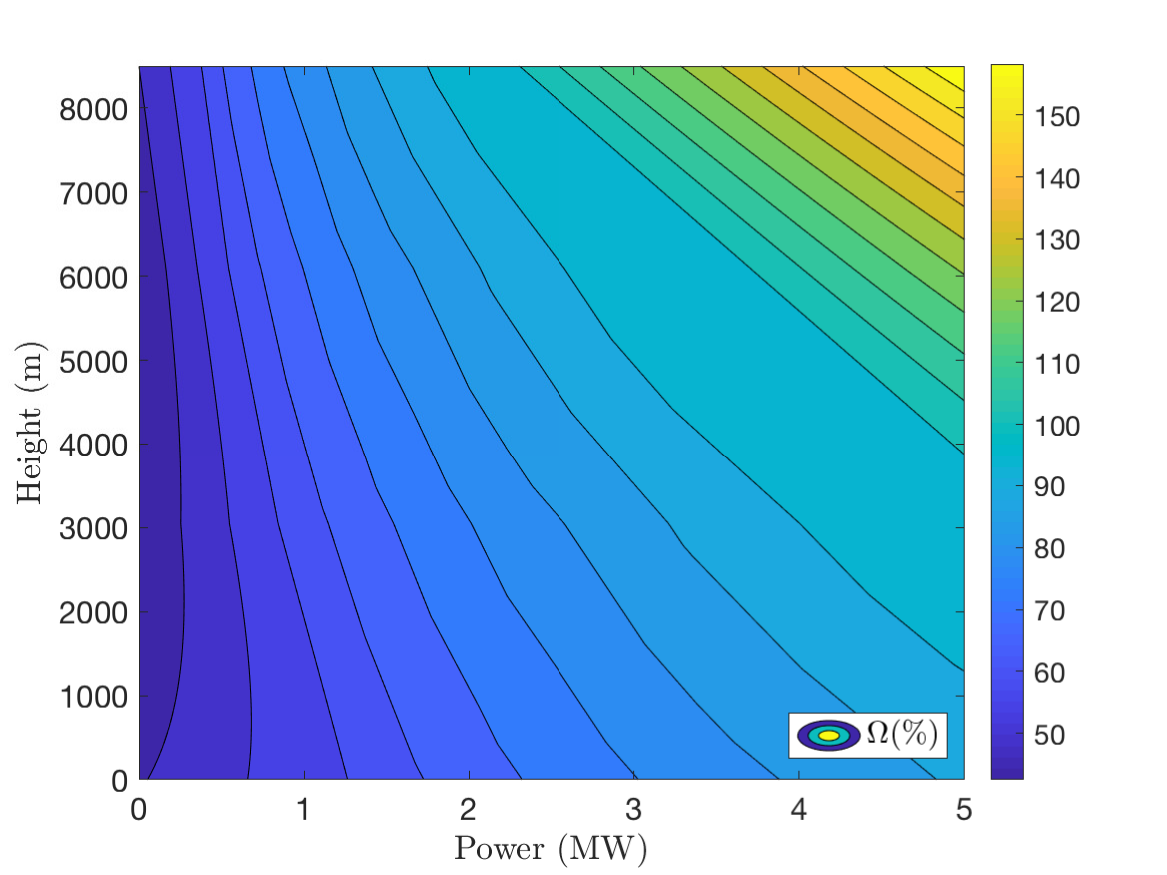}}
\caption{Contour plot relating drive power, altitude and non-dimensional rotation speed for a Mach number of 0.55.}
\label{fig:fanmap}
\end{figure}

The electric loss map coefficients $\kappa_{2,i},\kappa_{1,i},\kappa_{0,i}$ can be estimated $\forall i$ from these profiles. First, the  drive power $P_{\drv}$ is approximated a priori, e.g.~by assuming a conventional gas-turbine-powered flight. Then, the fan shaft rotation speed, $\omega_{\drv,i}$ (equal to the electric motor shaft rotation speed in both configurations), is interpolated from a precomputed look-up table relating measured shaft rotation speed, altitude and drive power at a given Mach number. For example, a Mach number of $0.55$ ($190$\,$\si{m/s}$ TAS) gives the relationship shown in Fig.~\ref{fig:fanmap}, which was obtained by scaling a proprietary fan design for the thrust range of the BAe 146 aircraft. The non-dimensional rotation speed $\Omega$ is thus estimated at a given altitude, Mach number and drive power using the map in Fig.~\ref{fig:fanmap}, and the shaft rotation speed is inferred from 
\[\omega_{drv} = \smash{\frac{156.7}{100}\frac{\pi}{30} \Omega \sqrt{T_{in}}},\]
where $T_{in}=T_{0}(h) + v^2/2 c_p$ is the temperature at inlet of the fan, $c_p=1005$ $\si{J K^{-1} kg^{-1}}$ is the specific heat of air at constant pressure and $T_0(h)$ is the temperature of air at altitude $h$. Finally, the coefficients are interpolated from a precomputed record of losses in the electric motor as a function of $\omega_\drv$.


\mdseditbis{The gas turbine fuel map and generator loss map used in this study are approximately linear ($\beta_{2} \approx 0$ and $\nu_{2} \approx 0$) for the range of power conditions considered, and the coefficients are constant as discussed in Section~\ref{sec:convex}.}

\subsection{Results}
\label{sec:results}
The mission is simulated in both configurations with sampling interval $\delta=60$ $\si{s}$ over a one-hour shrinking horizon by solving the optimisation problem \eqref{eq:min2} at each time step and implementing the first element of the optimal power split sequence as an MPC law. The tolerance is set to $\epsilon_{\text{rel}}=5\mathrm{e}{-6}$, $\epsilon_{\text{abs}}=0$ and the penalty parameters are updated at intervals of $F_\sigma=500$ iterations. 

\begin{figure}
    \centering
    \includegraphics[width = 0.4\textwidth,trim=4cm 0cm 4cm 0cm]{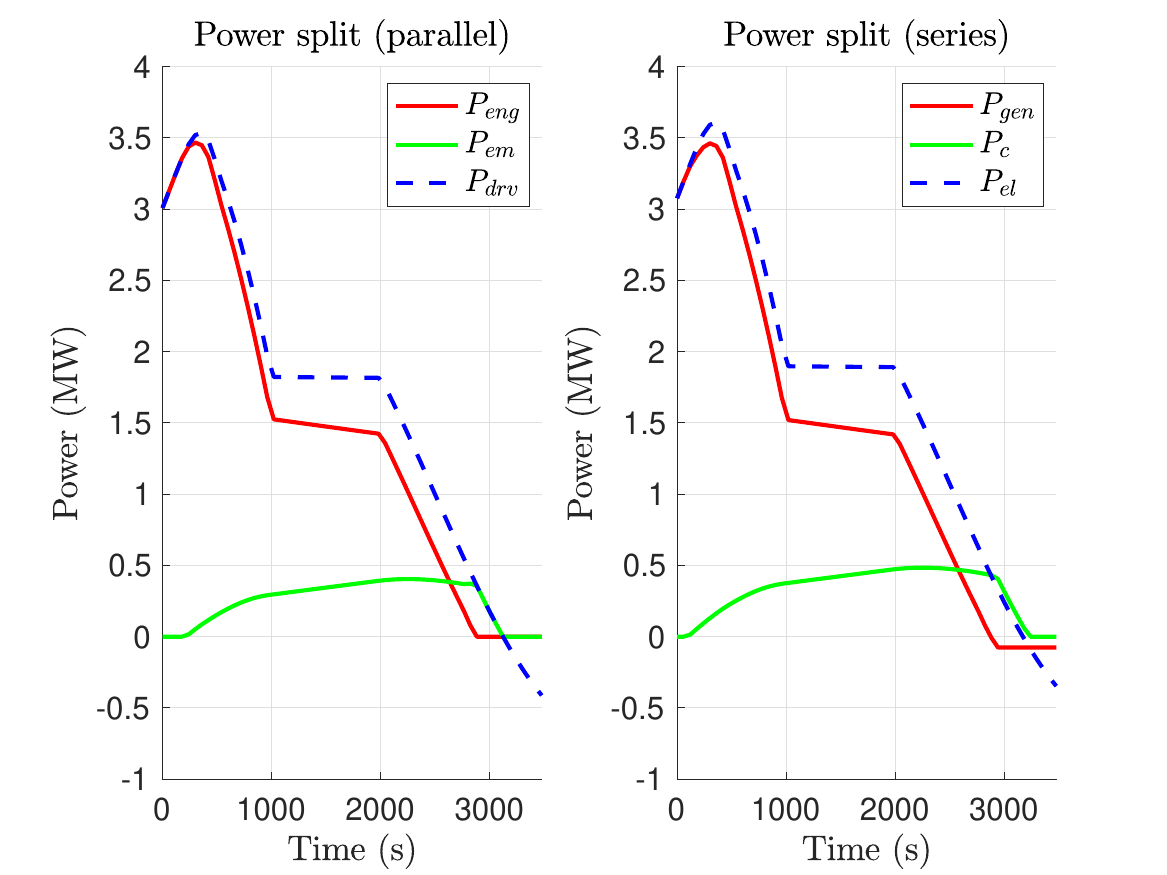}
    \caption{\mdseditbis{Closed-loop ADMM solution to the energy management problem in the parallel and series configurations, shown for 1 system (4 overall).}}
    \label{fig:parallel_series}
\end{figure}

\begin{figure}
    \centering
    \includegraphics[width = .4\textwidth,trim=4cm 0cm 4cm 1cm]{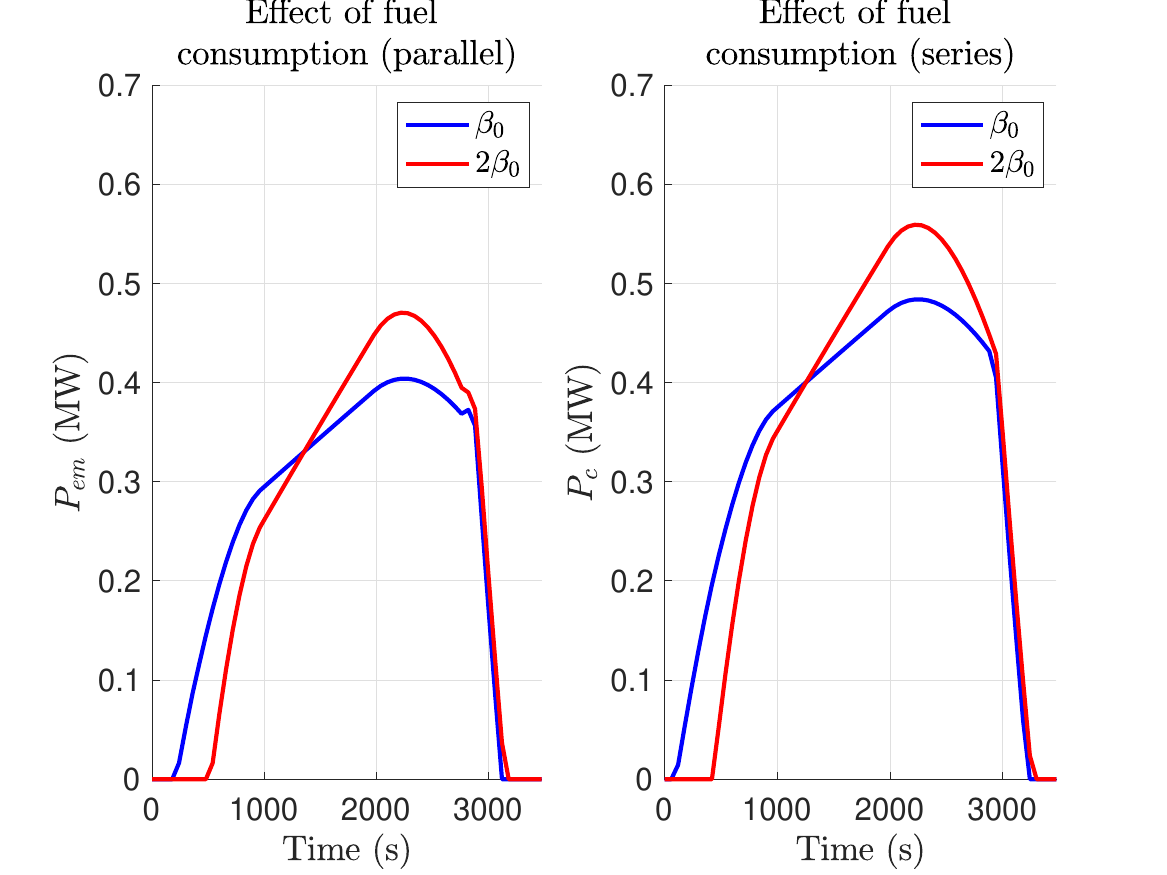}
    \caption{\mdsedittris{Effect of changing fuel map coefficient $\beta_0$ (single system).}}
    \label{fig:fuel}
\end{figure}

\mdseditbis{The closed-loop ADMM solution to the energy management control strategy is shown in Fig.~\ref{fig:parallel_series} for both parallel and series configurations. The solutions are for a single propulsion system  (so all quantities should be multiplied by $n=4$ to obtain the results for the whole aircraft). The plots represent the evolution of the relevant power terms in the power balance equations \eqref{eq:P_parallel} and \eqref{eq:P_series}: $P_\drv$, $P_\el$, $P_\gt$ for the parallel configuration and $P_\e$, $P_\c$, $P_\gen$ for the corresponding terms in the series configuration. It should be noted that the solution is similar in both configurations.} A striking feature of the solutions is the tendency to allocate more electrical power at the end of the flight. 
An intuitive explanation for this phenomenon is that the fuel burnt by using the gas turbine at the beginning of the flight reduces the mass of the aircraft, consequently reducing the power required to be produced by the fan later in the flight. This effect is amplified if the rate of fuel consumption is increased, as seen in Fig.~\ref{fig:fuel} comparing the electrical power profiles with different fuel consumption coefficients ($\beta$). 

\begin{figure}
    \centering
    \includegraphics[width = .4\textwidth,trim=4cm 0cm 4cm 0cm]{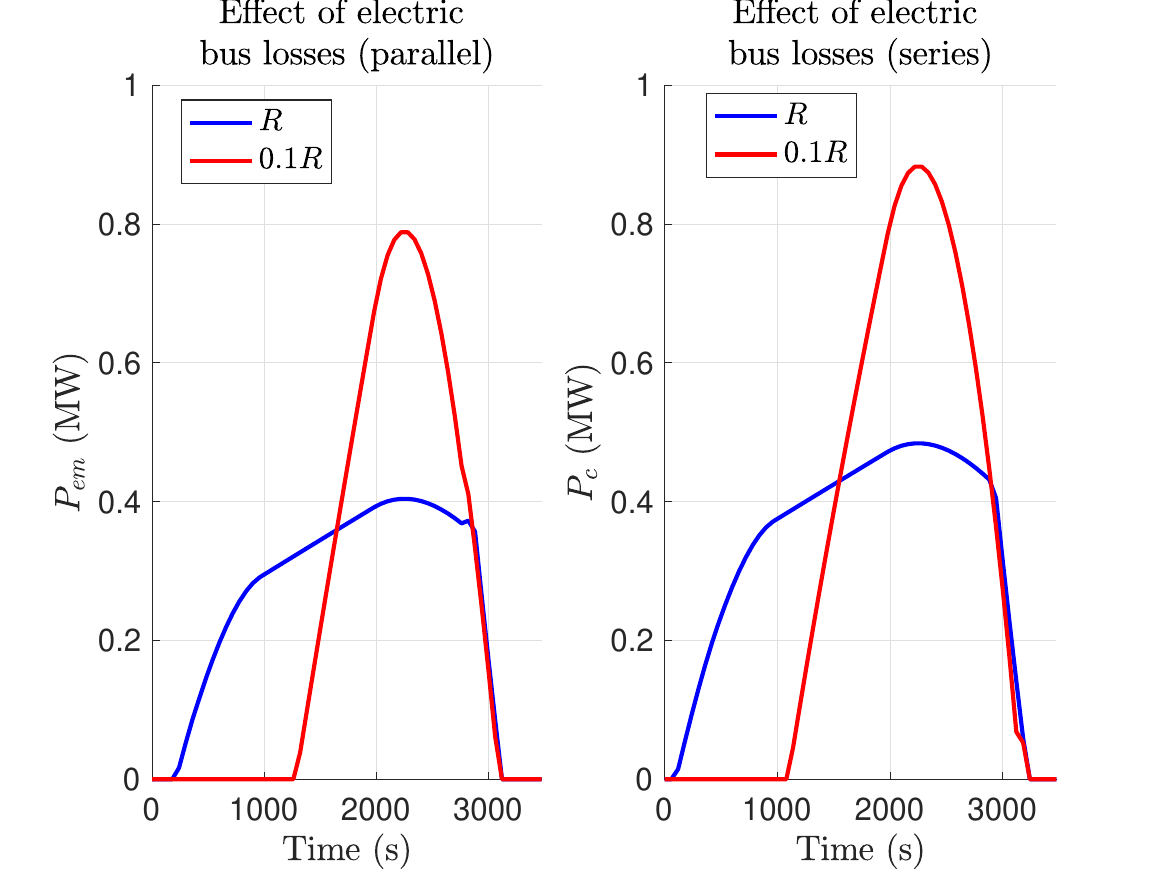}
    \caption{\mdseditbis{Effect of changing battery resistance \(R\) (single system)}.}
    \label{fig:losses}
\end{figure}

It should be noted that a concurrent effect arises from the losses in the battery electric bus. The nonlinear loss map $g$ between the battery chemical power $P_\b$ and effective power $P_\c$ tends to penalise large electrical power peaks thus flattening the electrical power distribution. This is seen in Fig.~\ref{fig:losses} which shows the electrical power profiles with different values for the battery equivalent circuit resistance: for smaller resistances the electrical losses at high power outputs is reduced so the power profile shows greater variation over time. 

\begin{figure}
    \centering
    \includegraphics[width = .4\textwidth,trim=4cm 0cm 4cm 0cm]{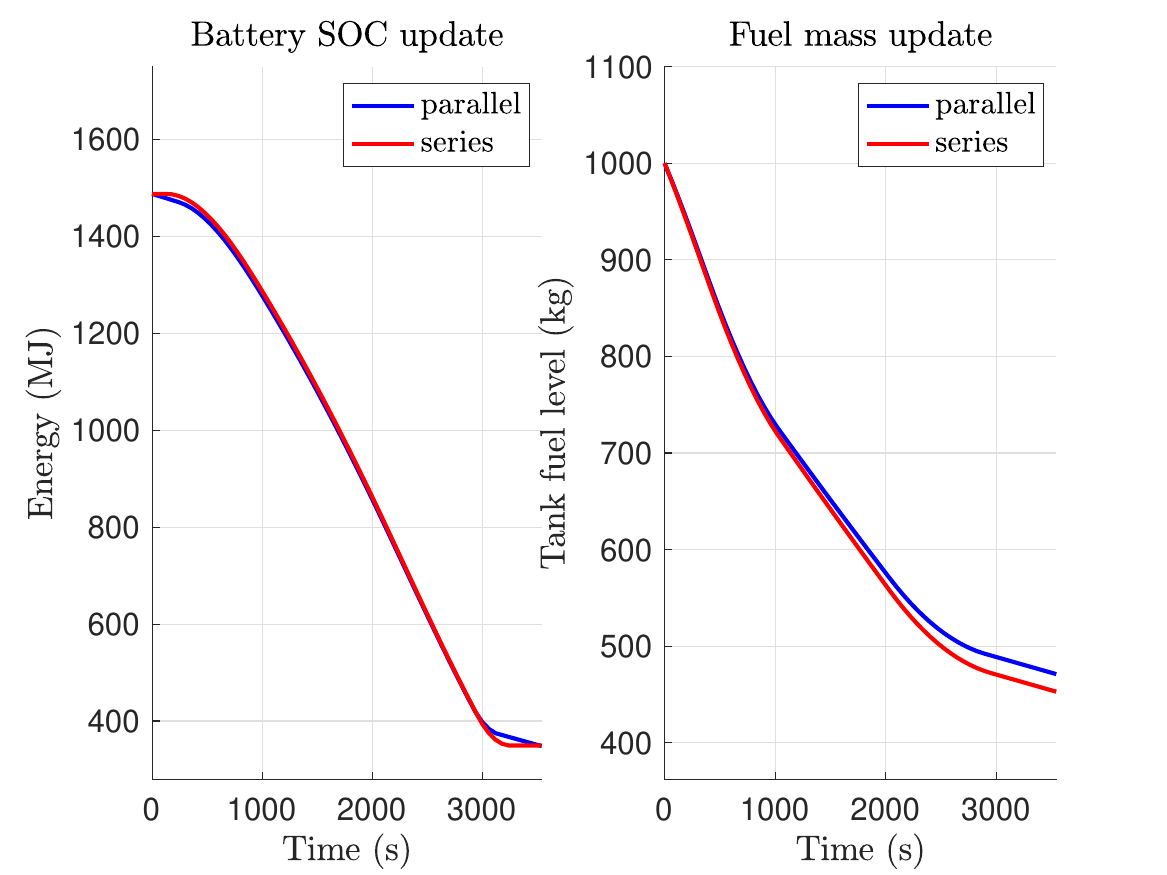}
    \caption{\mdseditbis{Comparison of battery and fuel consumption in the parallel and series configurations, shown for a single system.}}
    \label{fig:parallel_series2}
\end{figure}

Fig.~\ref{fig:parallel_series2} compares the evolution of battery SOC and fuel consumption for both configurations. As illustrated, the series propulsion architecture consumes slightly more fuel because it implements one more electric machine with associated losses, and so the electrical power is larger for a given flight profile. \finaledit{The selection of a particular configuration is thus motivated by a trade-off between efficiency and complexity of aero-mechanical integration.}

\begin{figure}
    \centering
    \includegraphics[width = .6\textwidth]{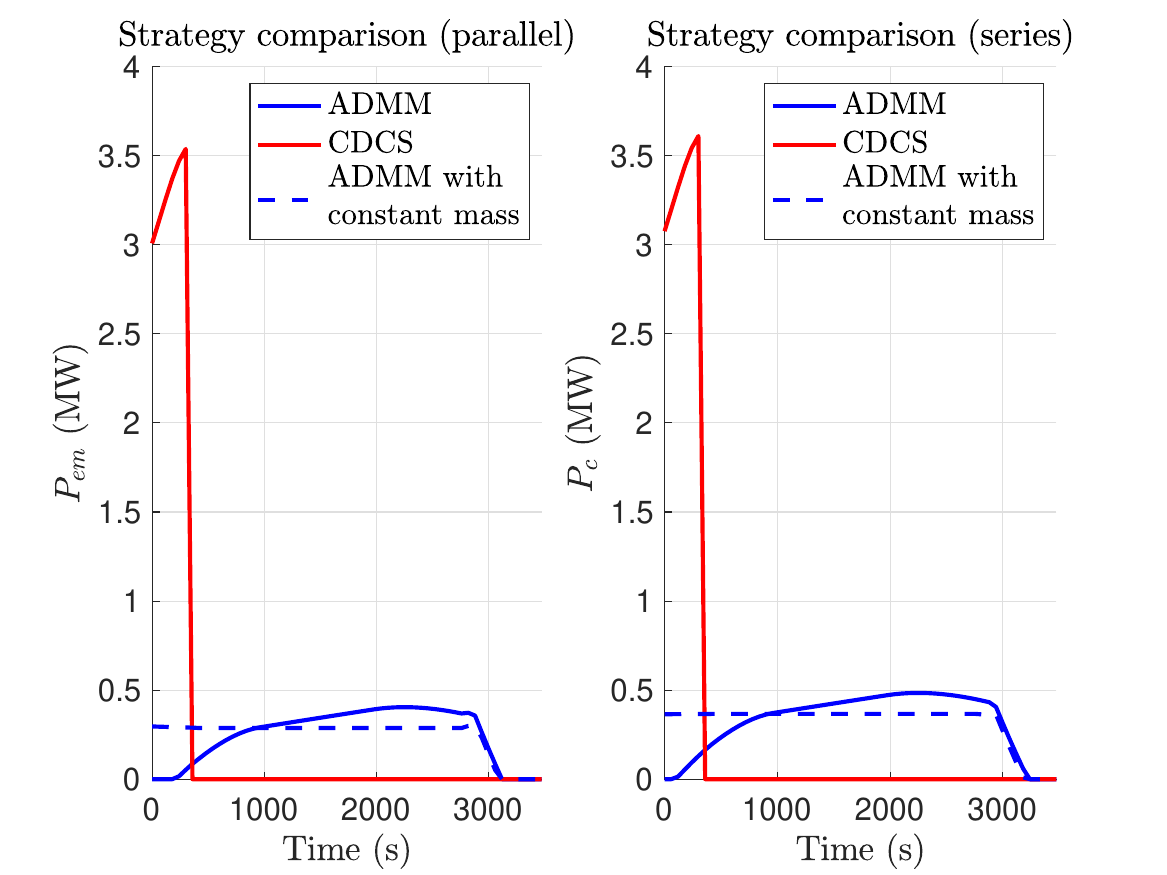}
    \caption{\mdseditbis{Comparison of ADMM, CDCS and ADMM with constant mass (parallel and series architectures, single system).}}
    \label{fig:strategies}
\end{figure}

The distribution of the electrical power over time is also illustrated in Fig.~\ref{fig:strategies} in comparison with other energy management strategies. The charge-depleting-charge-sustaining (CDCS) strategy is a heuristic that uses all the electrical energy at the beginning of the flight until the battery is depleted and then relies solely on the gas turbine for the remainder of the flight. Interestingly, the proposed ADMM-based approach is the antithesis of this strategy, allocating a non-negligible part of the electrical power at the end of the flight. The third strategy illustrated in Fig.~\ref{fig:strategies} uses the ADMM algorithm but ignores the aircraft mass variation. 
Interestingly, this (necessarily suboptimal) solution distributes the electrical power uniformly over the duration of the flight. In this case the strategy is dominated by the need to reduce electrical losses;
neglecting the aircraft mass variation means that the potential savings due to fuel burn early in the flight are not exploited.

The superiority of the presented variable mass ADMM solver over other strategies is shown in Table \ref{tab:fuel}. The heuristic CDCS strategy is used as a benchmark case. It is shown that the fuel savings with the mass-varying ADMM solver are superior to other strategies, in both parallel and series configuration. 
In the parallel configuration, the proposed energy management strategy achieves a fuel consumption of $2115$ $\si{kg}$, namely $1.7\%$ less than with CDCS.  \mdseditbis{Likewise, in the series configuration, a fuel consumption of $2188$ $\si{kg}$ is reported, which corresponds to a $1.9\%$ decrease over CDCS. }

As expected, the series architecture consumes more fuel than the parallel architecture. This is because series propulsion architectures employ an additional electrical machine and thus consume more electrical power due to inherent losses. Despite being less efficient, the series architecture has potential advantages in multi-propulsor configurations and is less mechanically complex than the parallel configuration. 

It should be noted that the same aircraft equipped with a conventional gas turbine propulsion system would burn $F_\gt = 2403$ $\si{kg}$ over the same scenario flight using the same models of powertrain components. However, in practice the conventional powertrain would be lighter since aviation fuels have a much higher energy density than batteries, so a direct comparison of fuel consumption is not possible.

\begin{table}
\centering
\begin{tabular}{lll|ll} 
\hline
                & \multicolumn{4}{c}{\textbf{Configuration}}                                                                                                \\ 
\cline{2-5}
                & \multicolumn{2}{c|}{\textbf{Parallel}}    & \multicolumn{2}{c}{\textbf{Series}}                                                           \\ 
\hline
\textbf{Method} & \textbf{Fuel (kg)} & \textbf{Saving (\%)} & \textbf{Fuel (kg)} & \textbf{Saving (\%)}                                                     \\ 
\hline
CDCS            & $2152$               & $-$                   & $2231$             & $-$                                                                      \\ 
\hline
Constant mass   & $2123$               & $1.3$                 & $2192$               & $1.7$                                                                      \\ 
\hline
Variable mass   & $2115$               & $1.7$                 & $2188$               & $1.9$                                                                     \\
\hline
\end{tabular}
\caption{\mdseditbis{Fuel comparison}}\label{tab:fuel}
\end{table}

Finally, we consider extensions of this case study to demonstrate the full potential of the proposed solver. We consider the same flight scenario but now assume that 1) the maximum gas turbine power is $\overline{P}_\gt=3$ MW, and 2) the propulsion unit is capable of converting negative drive power during descent into electricity to recharge the battery (i.e.~``windmilling''). This operation mode can be enforced by assuming a recovery efficiency $\eta_w$ and setting $\underline{P}_\b=\overline{P}_\b=g\bigl(h(\eta_w P_{\drv,i} )\bigr)$ for all time steps $i$ such that $P_{\drv,i} < 0$.

\begin{figure}
    \centering
    \includegraphics[width = 0.7\textwidth,trim=8cm 2cm 8cm 1cm]{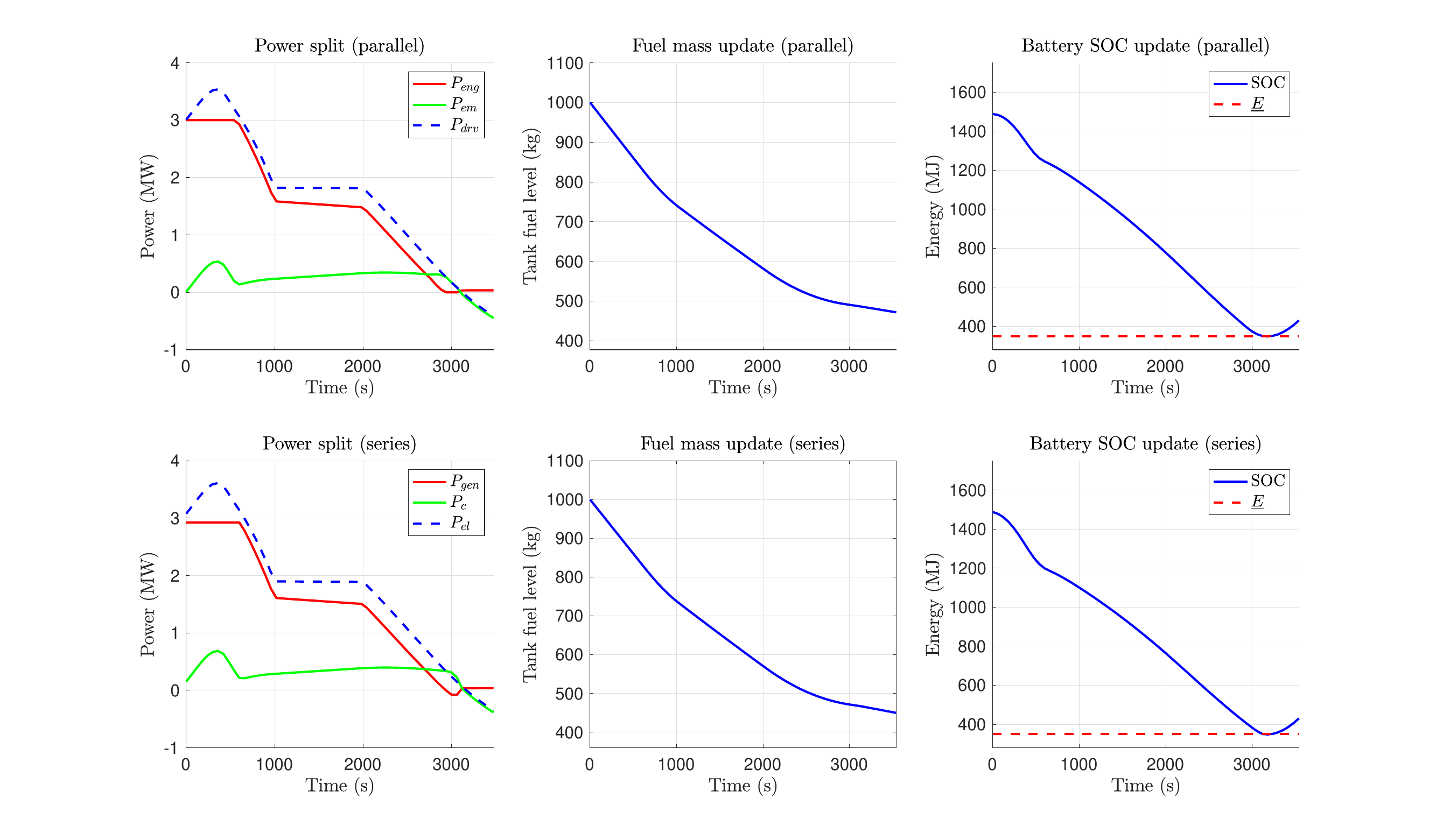}
    \caption{\mdsedit{Effect of windmilling and gas turbine saturation (single system)}.}
    \label{fig:saturation}
\end{figure}

Figure \ref{fig:saturation} shows the impact of these modifications. Gas turbine saturation causes more electrical energy to be allocated to the point at which  the gas turbine saturates. The potential for energy recovery via a windmilling mode is apparent at the end of the flight, where the drive power is negative, and the battery SOC increases during this part of the descent. It has been assumed for simplicity that the recovery process is ideal, that the electric motor can be operated as a generator and that the fan can be operated in reverse (requiring e.g.~a variable pitch fan). In practice, we would expect recovery efficiencies between $10\% - 20\%$ with current technology.

\subsection{Solver performance}

We next consider the convergence and robustness properties of the proposed ADMM solver. \finaledit{Instead of solving the optimisation problem at successive time steps to derive the MPC law (as in Section~\ref{sec:results}), we consider solving only one instance of the optimisation problem in order to simplify the analysis}. We show that the proposed solver is robust to changes in the flight profile, aircraft parameters and problem dimension.
The parallel configuration is considered here, all results being qualitatively equivalent for the series configuration. 

\begin{figure}
    \centering
    \includegraphics[width = .7\textwidth]{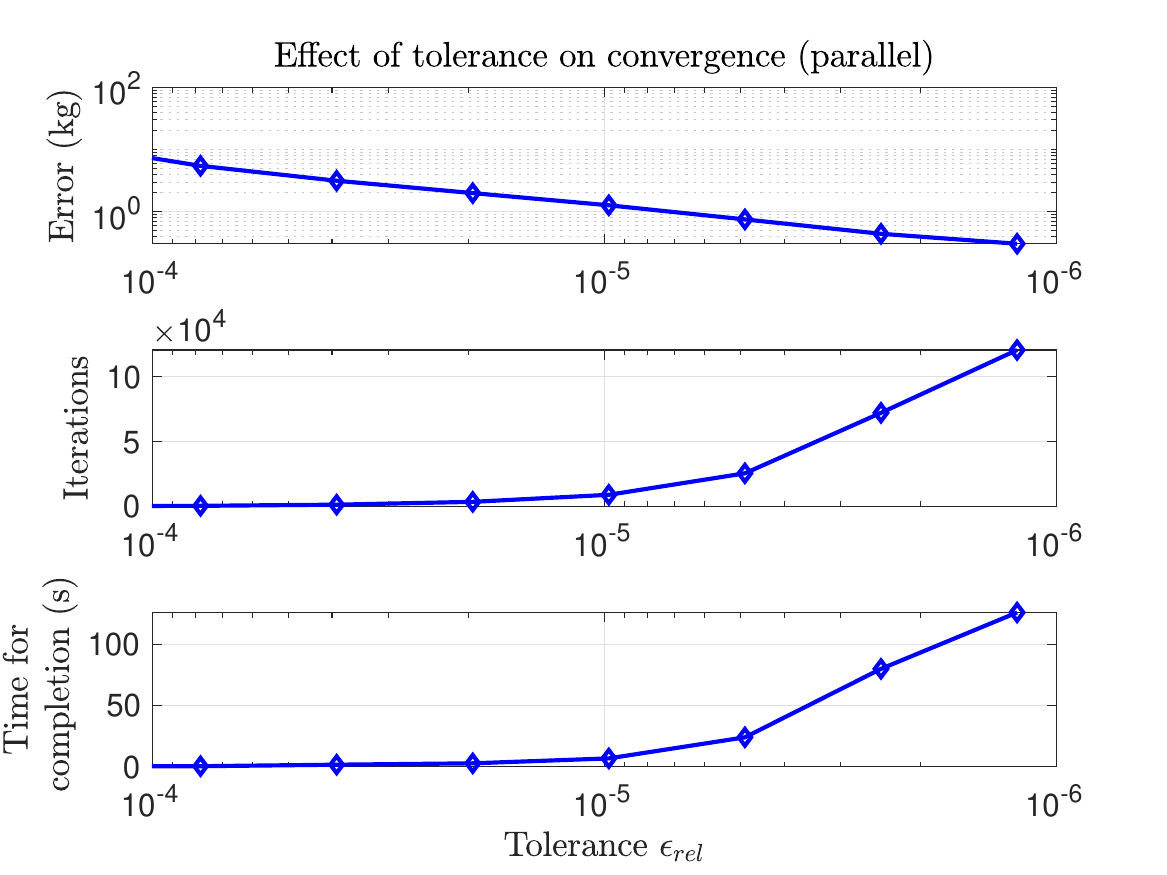}
    \caption{Effect of varying the relative tolerance on ADMM convergence.}
    \label{fig:convergence1}
\end{figure}

Section~\ref{sec:results} assumed fixed values for  parameters that influence convergence rate (tolerances, sample rate and penalty parameter update frequency). To investigate the effect of changing tolerances,
Fig.~\ref{fig:convergence1} shows accuracy relative to the optimal solution (obtained by solving the problem with optimisation package CVX and solver SDPT3 and comparing  total fuel consumption), number of iterations to completion, and computation time as a function of relative tolerance $\epsilon_\text{rel}$. The latter was varied while keeping other parameters constant (with $\epsilon_\text{abs}=0$, $F_\sigma = 10 ^5$, $\delta = 10 s$). As expected, as tolerance decreases, accuracy increases at the expense of a larger number of iterations and a consequent increase in computation. 

\begin{figure}
    \centering
    \includegraphics[width = .7\textwidth]{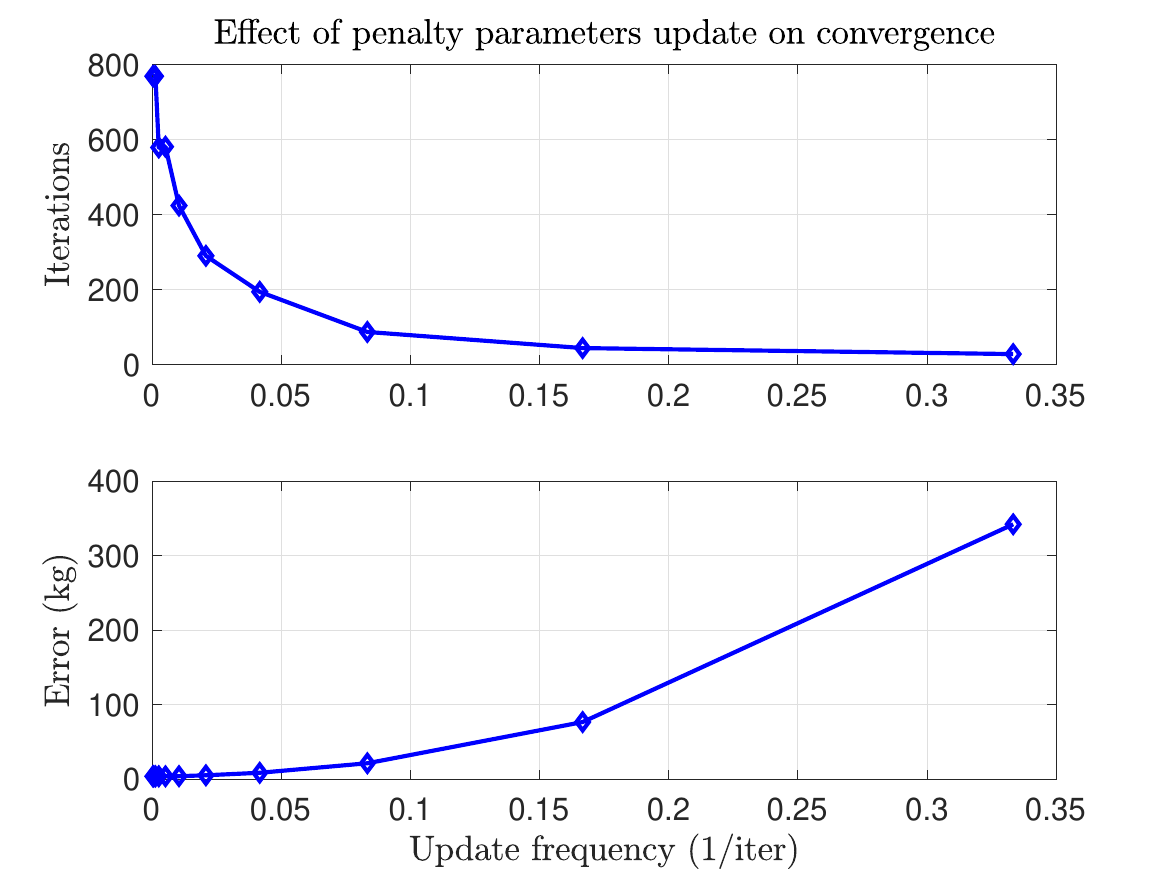}
    \caption{Effect of varying the penalty parameters update frequency on ADMM convergence.}
    \label{fig:convergence2}
\end{figure}

It is possible to reduce the tolerance without incurring additional computational cost if the ADMM algorithm is augmented with a penalty parameter update scheme as introduced in section \ref{sec:admm}. This is illustrated in Fig.~\ref{fig:convergence2}, which was obtained by varying the update frequency $1/F_\sigma$ while keeping other parameters constant (with $\epsilon_\text{abs}=0$, $\epsilon_\text{rel}=5 \times 10 ^{-5} $, $\delta = 10 s$). The number of iterations required (and consequently the computation time) decreases as the update frequency increases. However, this tends to decrease accuracy with respect to the CVX solution. The frequency update should thus be selected with care so as not to affect accuracy. 

\begin{figure}
    \centering
    \includegraphics[width = .7\textwidth]{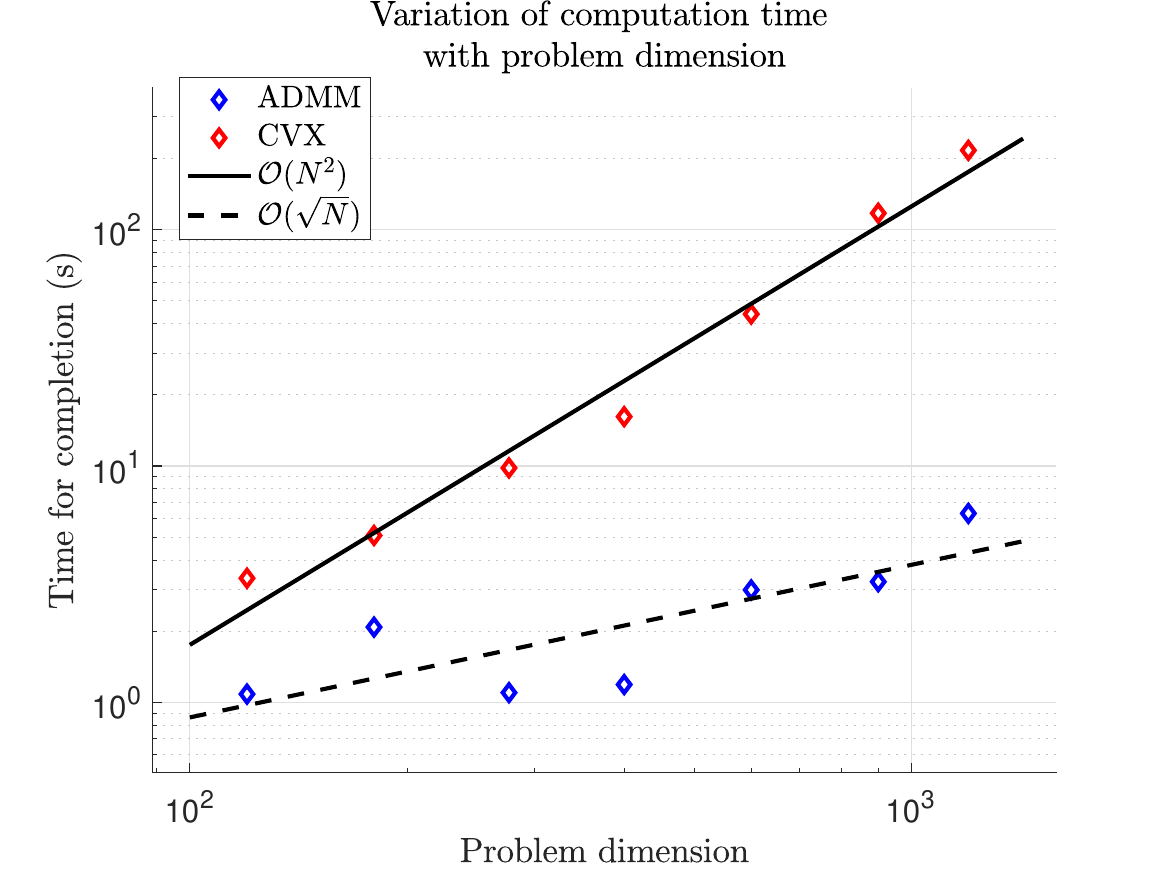}
    \caption{Effect of problem dimension on computation time.}
    \label{fig:dim}
\end{figure}

The influence of the sampling interval $\delta$ on computation time is shown in Fig.~\ref{fig:dim} by varying the problem dimension ($N = T / \delta$), with all other parameters kept constant ($\epsilon_\text{abs}=0$, $\epsilon_\text{rel}=5 \times 10 ^{-5} $, $F_\sigma = 50$). Computation time increases as problem dimension increases, however, the \finaledit{empirically observed} dependence is $\mathcal{O}(N^2)$ for CVX and $\mathcal{O}(\sqrt{N})$ for ADMM. Therefore the proposed solver provides significant computation time reduction relative to CVX, allowing longer prediction horizons and better real-time convergence. 

Experiments were performed to compare the proposed ADMM algorithm for the convex problem (\ref{eq:min2}) with direct solution of the nonconvex problem~(\ref{eq:min}). Retaining only assumption 1) from section \ref{sec:convex}, the nonconvex problem was solved using a general purpose nonlinear programming solver (fmincon \cite{MatlabOTB}) with $\delta = 60$\,s, which converged within $82$\,s. Under the same conditions ADMM (implemented in Matlab) converged within $0.5$\,s. To compare fmincon and ADMM solutions, a Monte Carlo simulation was conducted by solving $100$ problem instances with battery size randomly sampled from a uniform distribution. For each scenario the mean absolute error between the solutions ($P_b$) of both solvers was computed. The variance of the error distribution is $9.3 \times 10 ^{-5}$\, $\si{MW^2}$, showing good agreement between fmincon and ADMM.

\begin{figure}
    \centering
    \includegraphics[width = .7\textwidth]{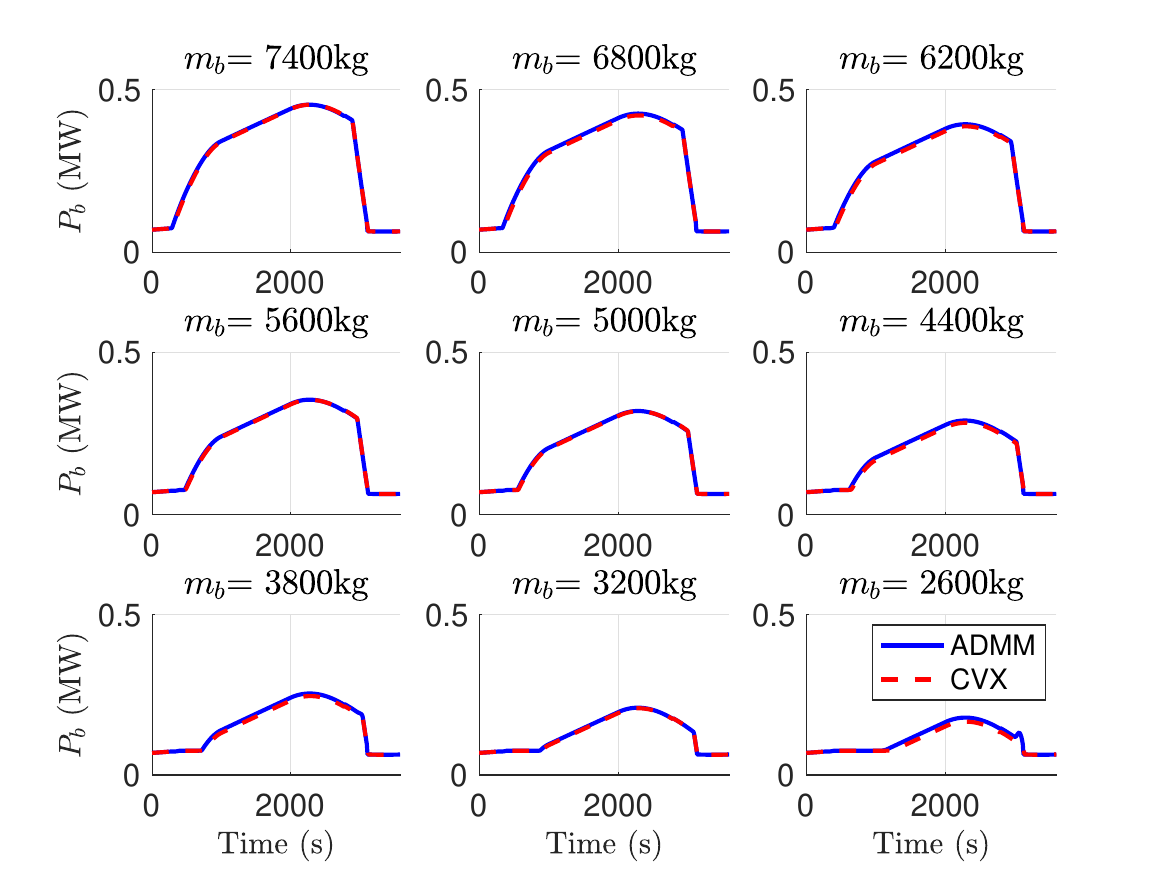}
    \caption{ADMM and CVX solutions for various battery masses.}
    \label{fig:robustness1}
\end{figure}

\begin{figure}
    \centering
    \includegraphics[width = .7\textwidth]{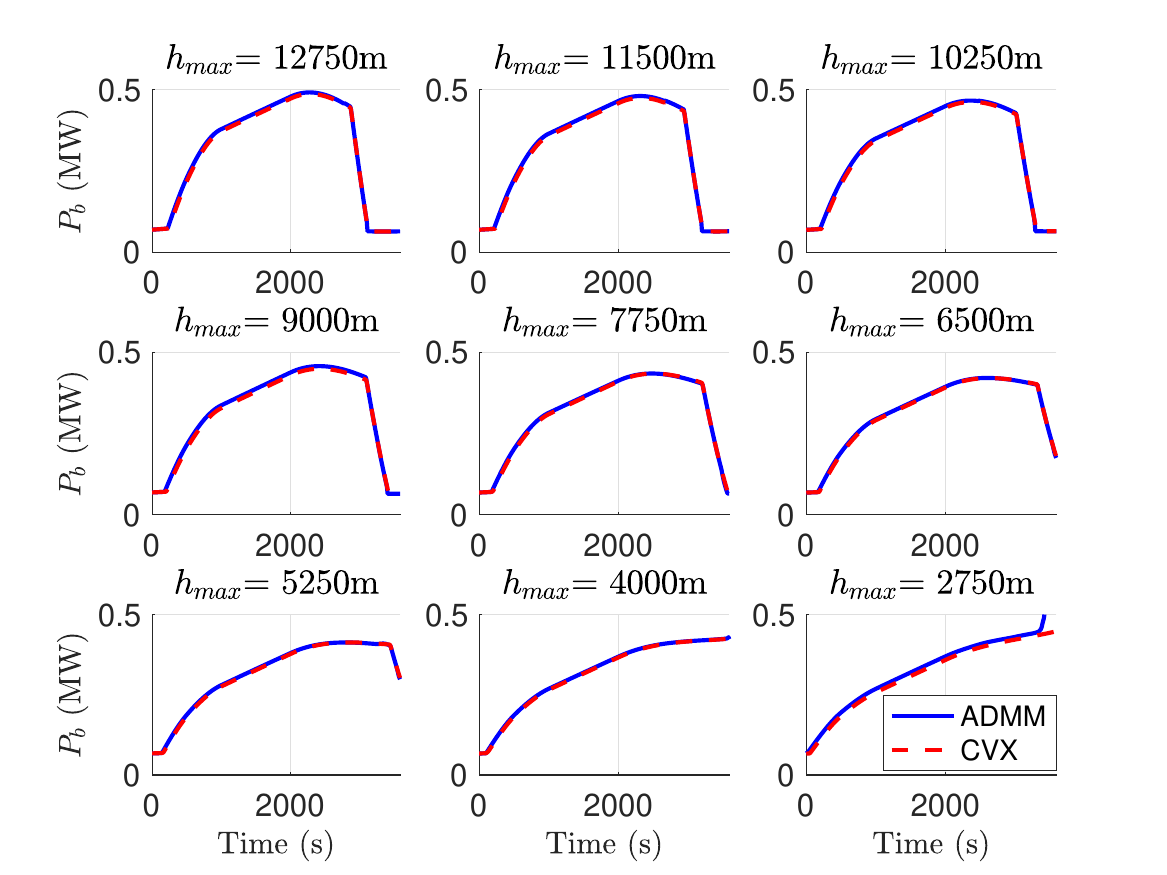}
    \caption{ADMM and CVX solutions for various maximum altitudes.}
    \label{fig:robustness2}
\end{figure}

\begin{figure}
    \centering
    \includegraphics[width = .7\textwidth]{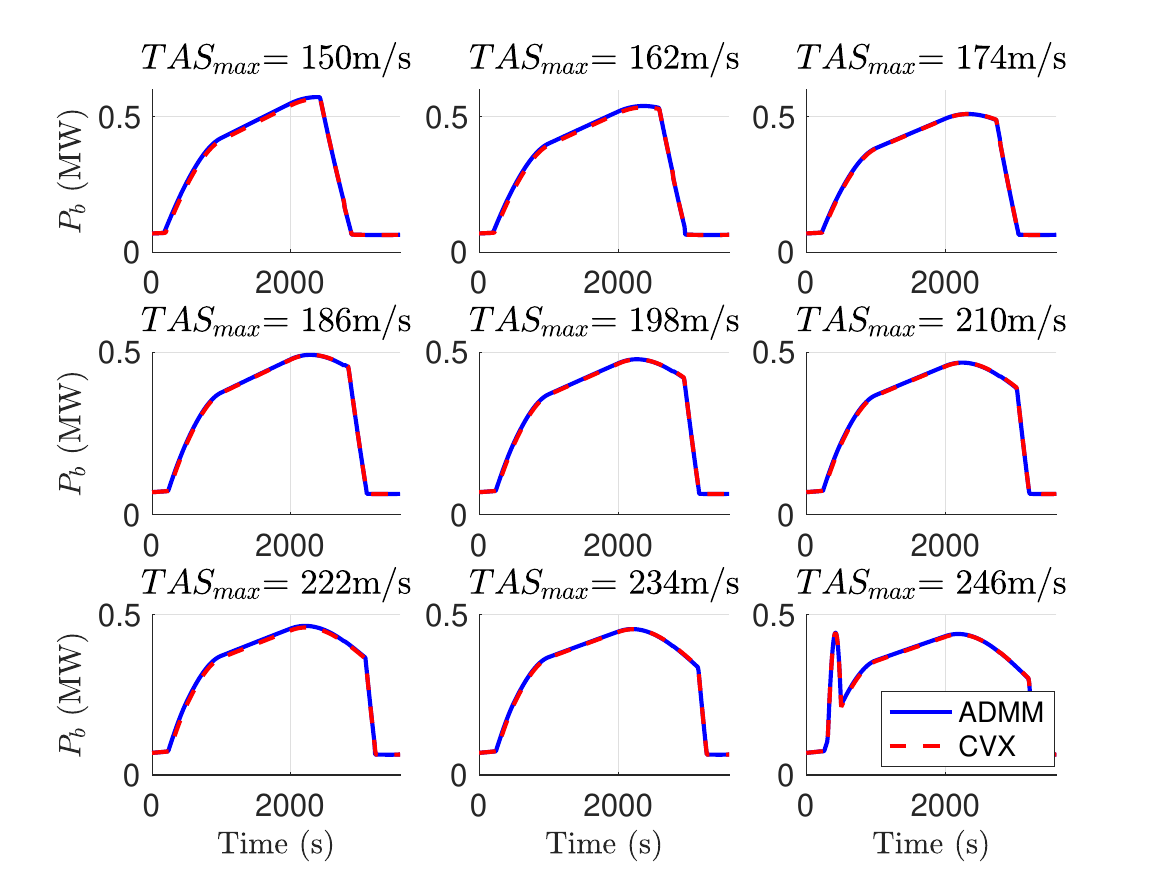}
    \caption{ADMM and CVX solutions for various maximum TAS.}
    \label{fig:robustness3}
\end{figure}

Finally, robustness to changes in the mission parameters is investigated in Figures \ref{fig:robustness1}-\ref{fig:robustness3} where CVX and ADMM solutions are compared for modified simulation scenarios and fixed convergence parameters ($\epsilon_\text{abs}=0$, $\epsilon_\text{rel}=5 \times 10 ^{-5} $, $F_\sigma = 50$, $\delta = 10 s$). These results show that the ADMM solution matches the solution obtained using CVX, thus demonstrating robustness to changes in problem-specific parameters.

\section{Conclusions}
\label{sec:conclusion}

This paper presents a fast and robust ADMM algorithm to solve the energy management problem for a hybrid electric aircraft in parallel and series configurations. A convex program is derived from the associated optimisation problem and the high degree of separability in the optimisation variables is exploited in the design of the solver. The ADMM solver was shown to produce similar results to the general purpose convex optimisation package CVX (with solver SDPT3) for a wide range of scenarios, while significantly outperforming CVX in terms of computation times. 
Significant fuel savings were achieved by comparison to  heuristic strategies. 

\mdseditbis{An extension of the proposed approach could optimise gas turbine speed given an estimate of its power output (possibly within an iterative scheme), removing the need for the assumption on gas turbine speed for the series configuration}.
\mdsedittris{Another extension would be to investigate robustness of the proposed approach to power demand disturbances.  Although the flight path is fixed and the aircraft flight dynamics are prescribed in the MPC optimisation, the predicted power demand is likely to be inexact and this would introduce disturbance terms into the dynamics of battery SOC and fuel mass.} 
Future work will also investigate the application of the proposed algorithm to solve the energy management problem for other types of hybrid vehicles (e.g. hybrid VTOL aircraft with applications to urban air mobility). \mdsedittris{Finally, the principles developed here for energy management could be applied to the problem of optimal design and sizing of powertrain components.}

\bibliography{biblio} 
\bibliographystyle{plain}

\end{document}